\documentclass[11pt]{article}
\usepackage{mathrsfs}
\usepackage{amsthm}
\usepackage{amssymb}
\usepackage{amsmath}
\usepackage{graphicx}
\usepackage{color}
\usepackage{amsfonts}
\usepackage{float}
\usepackage{cite}
\usepackage [latin1]{inputenc}
\usepackage[text={140mm,210mm},left=35mm,vmarginratio=1:1]{geometry}
\newtheorem{theorem}{Theorem}[section]

\newtheorem{lemma}[theorem]{Lemma}

\newtheorem{corollary}[theorem]{Corollary}

\numberwithin{equation}{section}
\normalsize

\begin{document}
\title{\textbf{Hydrodynamics of a class of $N$-urn linear systems}}

\author{Xiaofeng Xue \thanks{\textbf{E-mail}: xfxue@bjtu.edu.cn \textbf{Address}: School of Science, Beijing Jiaotong University, Beijing 100044, China.}\\ Beijing Jiaotong University}

\date{}
\maketitle

\noindent {\bf Abstract:} In this paper we are concerned with hydrodynamics of a class of $N$-urn linear systems, which include voter models, pair-symmetric exclusion processes and binary contact path processes on $N$ urns as special cases. We show that the hydrodynamic limit of our process is driven by a $\left(C[0,1]\right)^\prime$-valued linear ordinary differential equation and the fluctuation of our process, i.e, central limit theorem from the hydrodynamic limit, is driven by a $\left(C[0, 1]\right)^\prime$-valued Ornstein-Uhlenbeck process. To derive above main results, we need several replacement lemmas.
An extension in linear systems of Chapman-Kolmogorov equation plays key role in proofs of these replacement lemmas.

\quad


\section{Introduction}\label{section one}
In this paper, we are concerned with a class of $N$-urn linear systems. The so-called $N$-urn linear system is a continuous-time Markov process $\{X_t\}_{t\geq 0}$ with state space $\mathbb{X}\subseteq [0,+\infty)^N$, where $N\geq 1$ is a given integer. For any $x\in [0, +\infty)^{N}$ and $1\leq i\leq N$, let $x(i)$ be the $i$th coordinate of $x$, then the transition rates of $\{X_t\}_{t\geq 0}$ is given as follows. For any $t\geq 0$ and $1\leq i\neq j\leq N$,
\begin{equation}\label{equ 1.1 transition rates}
X_t\rightarrow
\begin{cases}
X_t^{i, +} & \text{~at rate~} b\left(\frac{i}{N}\right),\\
X_t^{(i,j)} & \text{~at rate~} \frac{1}{N}\lambda\left(\frac{i}{N}, \frac{j}{N}\right),
\end{cases}
\end{equation}
where
\[
x^{i, +}(k)=
\begin{cases}
x(k) & \text{~if~}k\neq i,\\
c\left(\frac{i}{N}\right)x(i) & \text{~if~} k=i
\end{cases}
\]
and
\[
x^{(i,j)}(k)=
\begin{cases}
x(k) & \text{~if~}k\neq i, j,\\
a_1\left(\frac{i}{N}, \frac{j}{N}\right)x(i)+a_2\left(\frac{i}{N}, \frac{j}{N}\right)x(j) & \text{~if~} k=i,\\
a_3\left(\frac{i}{N}, \frac{j}{N}\right)x(i)+a_4\left(\frac{i}{N}, \frac{j}{N}\right)x(j) & \text{~if~} k=j
\end{cases}
\]
for any $x\in [0, +\infty)^N$ and $1\leq i\neq j\leq N$, where $b,c\in C[0, 1]$ and $\lambda, a_1, a_2, a_3, a_4\in C([0, 1]\times [0, 1])$ are all non-negative. From now on, we write $X_t$ as $X_t^N$ when we need to point out the integer $N$. $\{X_t^N\}_{t\geq 0}$ can be defined equivalently via its generator. By \eqref{equ 1.1 transition rates}, the generator $\mathcal{L}_N$ of $\{X_t^N\}_{t\geq 0}$ is given by
\begin{equation}\label{equ 1.2 generator}
\mathcal{L}_Nf(x)=\sum_{i=1}^Nb\left(\frac{i}{N}\right)\left(f(x^{i,+})-f(x)\right)+\frac{1}{N}\sum_{i=1}^N\sum_{j\neq i}\lambda\left(\frac{i}{N}, \frac{j}{N}\right)\left(f(x^{(i,j)})-f(x)\right)
\end{equation}
for any $x\in [0, +\infty)^N$ and $f$ from $[0, +\infty)^N$ to $\mathbb{R}$. $\{X_t^N\}_{t\geq 0}$ is a special case of the linear system introduced in Chapter 9 of \cite{Lig1985} since $X_t^N$ jumps to each probable linear transformation of itself at a constant rate.

Our model includes several important examples when $\lambda, a_1, a_2, a_3, a_4, b, c$ are specially given.

\textbf{Example 1} \emph{$N$-urn voter model}. When $b=0, a_1=0, a_2=1, a_3=0, a_4=1$ and $X_0^N\in \{0, 1\}^N$, $\{X_t^N\}_{t\geq 0}$ reduces to the voter model (see Chapter 5 of \cite{Lig1985} or Part {\rm \uppercase\expandafter{\romannumeral2}} of \cite{Lig1999}) on $N$ urns, where each urn has an option $0$ or $1$ and the $i$th urn adopts the option of the $j$th one at rate $\frac{1}{N}\lambda\left(\frac{i}{N}, \frac{j}{N}\right)$.

\quad

\textbf{Example 2} \emph{Pair-symmetric $N$-urn exclusion process}. When $b=0, a_1=0, a_2=1, a_3=1, a_4=0$ and $X_0^N\in \{0, 1\}^N$, $\{X_t^N\}_{t\geq 0}$ reduces to a special case of the exclusion process (see Chapter 7 of \cite{Lig1985} or Part {\rm \uppercase\expandafter{\romannumeral3}} of \cite{Lig1999}) on $N$-urns, where each urn has at most one particle and a particle jumps from the $i$th (resp. $j$th) urn to the $j$th (resp. $i$th) one at rate
\[
\widetilde{\lambda}_N\left(\{i,j\}\right)=\frac{1}{N}\left(\lambda\left(\frac{i}{N}, \frac{j}{N}\right)+\lambda\left(\frac{j}{N}, \frac{i}{N}\right)\right)
\]
when the $j$th (resp. $i$th) urn is vacant. We call this model pair-symmetric exclusion process since the rate at which a particle jumps between a pair of urns does not depend on which one being the starting urn.

\quad

\textbf{Example 3} \emph{$N$-urn binary contact path process}. When $a_1=a_2=1, a_3=0, a_4=1, c=0$ and $X_0^N\in \{0,1,2,\ldots\}^N$, $\{X_t^N\}_{t\geq 0}$ reduces to the binary contact path process (BCPP) introduced in \cite{Gri1983} on $N$-urns,  where $X_t^N(i)$ is the seriousness of an infection disease on the $i$th urn. The $i$th urn becomes healthy at rate $b\left(\frac{i}{N}\right)$ and is further infected by the $j$th one at rate $\frac{1}{N}\lambda\left(\frac{i}{N}, \frac{j}{N}\right)$. When the $i$th box is infected by the $j$th one, the seriousness of the ill on the $i$th box is added with that of the $j$th one. The binary contact path process is an auxiliary tool for the investigation of the contact process (CP) (see Chapter 6 of \cite{Lig1985} or Part {\rm \uppercase\expandafter{\romannumeral1}} of \cite{Lig1999}) since the process $\{Y_t\}_{t\geq 0}$ with state space $\{0, 1\}^N$ defined as
\[
Y_t(i)=
\begin{cases}
1 & \text{~if~}X_t^N(i)>0,\\
0 & \text{~if~}X_t^N(i)=0
\end{cases}
\]
for each $1\leq i\leq N$ is a version of the contact process. According to the coupling relationship between BCPP and CP, several upper bounds of critical values of contact processes on $\mathbb{Z}^d$ with $d\geq 3$ are given. For mathematical details, see \cite{Gri1983} and \cite{Xue2018}.

\quad

In this paper, we are concerned with the hydrodynamic limit of $\{X_t^N\}_{t\geq 0}$ and the corresponding fluctuation, i.e, the law of large numbers of $\frac{1}{N}\sum_{i=1}^NX_t^N(i)\delta_{\frac{i}{N}}(du)$ and the central limit theorem from it, where $\delta_a(du)$ is the Dirac measure concentrated on $a$. The motivation for the investigation of the hydrodynamic limit is to derive an differential equation, which is usually a real-valued PDE or an ODE on a Banach space, describing the macroscopic behavior of the microscopic density field of the stochastic model. For a comprehensive reading of this topic, see \cite{kipnis+landim99}.

Investigations of hydrodynamics of linear systems date back to 1980s. In \cite{PresuttiSpohn83}, the hydrodynamic limit and corresponding fluctuation of voter models on lattices are discussed. In detail, it is shown in \cite{PresuttiSpohn83} that, after a diffusive scaling, the hydrodynamic limit of the voter model on $\mathbb{Z}^d$ with $d\geq 3$ is driven by a heat equation and the fluctuation from the above limit is driven by a generalized Ornstein-Uhlenbeck process introduced in \cite{Holley1978}, which is tempered distribution-valued. Recently, hydrodynamics of binary contact path processes on lattices are discussed. Similar with main theorems in \cite{PresuttiSpohn83}, it is shown in \cite{XueZhao2020} and \cite{XueZhao2021} that, after a diffusive scaling, the hydrodynamic limit of the BCPP on $\mathbb{Z}^d$ with infection rate $\lambda$ is driven by a heat equation while the corresponding fluctuation is driven by a tempered distribution-valued O-U process when $d$ and $\lambda$ are sufficiently large.

So far, investigations of hydrodynamics of interacting particle systems on $N$-urns are mainly concerned with the Ehrenfest model, where particles perform independent random walks between $N$-urns. Hydrodynamic limits and corresponding fluctuations, large and moderate deviations of generalized $N$-urn Ehrenfest models are given in \cite{Xue2022} and \cite{Ren2021}. It is more important to investigate cases where interactions are among particles. So in this paper we discuss the linear system on $N$-urns, which includes important models such as voter models, symmetric exclusion processes and binary contact path processes. It is also interesting to investigate hydrodynamics of other interacting particle systems on $N$-urns such as asymmetric exclusion process, zero-range processes, misanthrope processes, branching random walks. We will work on these models as further investigations.

\section{Main results}\label{section two}
In this section we give our main results. We first introduce some notations and definitions for later use. For $k\in\{1,2,\ldots\}$, we use $\left(C\left([0, 1]^k\right)\right)^\prime$ to denote the dual of $C\left([0,1]^k\right)$. For $k\in \{1,2,\ldots\}$ and any $f\in C\left([0, 1]^k\right)$, we use $\|f\|_\infty$ to denote the $l_\infty$-norm of $f$, i.e.,
\begin{equation}\label{equ functionInfinityNorm}
\|f\|_\infty=\sup_{u\in [0, 1]^k}|f(u)|.
\end{equation}
For $k\in \{1,2,\ldots\}$ and any $\nu\in \left(C\left([0, 1]^k\right)\right)^\prime$, we define
\begin{equation}\label{equ functionDualInifinityNorm}
\|\nu\|=\sup_{f\in C\left([0, 1]^k\right): ~\|f\|_\infty=1}|\nu(f)|.
\end{equation}
For $k_1, k_2\in \{1,2,\ldots\}$ and any bounded linear operator $\mathcal{A}$ from $C\left([0, 1]^{k_1}\right)$ to $C\left([0, 1]^{k_2}\right)$, we denote by $\mathcal{A}^{*}$ the linear operator from $\left(C\left([0, 1]^{k_2}\right)\right)^\prime$ to $\left(C\left([0, 1]^{k_1}\right)\right)^\prime$ such that
\[
(\mathcal{A}^{*}\nu)(f)=\nu(\mathcal{A}f)
\]
for any $f\in C\left([0, 1]^{k_1}\right)$ and $\nu\in \left(C\left([0, 1]^{k_2}\right)\right)^\prime$. For any $t\geq 0$, we use $\mu_t^N$ to denote $\frac{1}{N}\sum_{i=1}^NX_t^N(i)\delta_{\frac{i}{N}}(du)$, where $\delta_a(du)$ is the Dirac measure concentrated on $a$, then $\mu_t^N$ can be identified with a random element in $\left(C[0, 1]\right)^\prime$ such that
\[
\mu_t^N(f)=\frac{1}{N}\sum_{i=1}^NX_t^N(i)f\left(\frac{i}{N}\right)
\]
for any $f\in C[0, 1]$. Throughout this paper, we adopt the following initial assumption.

\textbf{Assumption} (A): $\{X_0^N(i)\}_{1\leq i\leq N}$ are independent and
\[
P\left(X_0^N(i)=1\right)=\phi\left(\frac{i}{N}\right)=1-P\left(X_0^N(i)=0\right)
\]
for each $1\leq i\leq N$, where $\phi\in C[0, 1]$.

We define $P_1$ as the linear operator from $C[0, 1]$ to $C[0, 1]$ such that
\begin{align*}
P_1f(u)=&b(u)\left(c(u)-1\right)f(u)+f(u)\int_0^1\lambda(u, v)\left(a_1(u, v)-1\right)dv\\
&+\int_0^1\lambda(u, v)a_3(u, v)f(v)dv+f(u)\int_0^1\lambda(v, u)\left(a_4(v, u)-1\right)dv\\
&+\int_0^1\lambda(v, u)f(v)a_2(v, u)dv
\end{align*}
for any $f\in C[0, 1]$ and $u\in [0, 1]$. According to the theory of ordinary differential equations on Banach spaces (see Chapter 19 of \cite{Lang}), we have the following lemma.

\begin{lemma}\label{lemma 2.1 ODE of Hydrodynamic}
There exists a unique solution $\{\mu_t\}_{t\geq 0}$ to the $\left(C[0, 1]\right)^\prime$-valued ordinary differential equation
\begin{equation}\label{equ definition of hydro ODE}
\begin{cases}
\frac{d}{dt}\mu_t=P_1^{*}\mu_t \text{~for~}t\geq 0,\\
\mu_0(du)=\phi(u)du
\end{cases}
\end{equation}
under the norm $\|\cdot\|$ given in \eqref{equ functionDualInifinityNorm}. Furthermore, $\mu_t(du)=\rho(t,u)du$ for all $t\geq 0$, where $\{\rho(t,\cdot)\}_{t\geq 0}$ is the unique solution to the $C[0,1]$-valued ordinary differential equation
\begin{equation}\label{equ definition of hydro density}
\begin{cases}
&\frac{d}{dt}\rho(t, \cdot)=b(\cdot)\left(c(\cdot)-1\right)\rho(t, \cdot)+\rho(t,\cdot)\int_0^1\lambda(v, \cdot)\left(a_4(v, \cdot)-1\right)dv\\
&\text{\quad\quad\quad\quad~}+\rho(t, \cdot)\int_0^1\lambda(\cdot, v)\left(a_1(\cdot, v)-1\right)dv+\int_0^1\lambda(\cdot, v)a_2(\cdot, v)\rho(t, v)dv\\
&\text{\quad\quad\quad\quad~}+\int_0^1\lambda(v, \cdot)a_3(v, \cdot)\rho(t, v)dv,\\
&\rho(0,\cdot)=\phi
\end{cases}
\end{equation}
under the norm $\|\cdot\|_\infty$ given in \eqref{equ functionInfinityNorm}.
\end{lemma}

We prove Lemma \ref{lemma 2.1 ODE of Hydrodynamic} in Appendix \ref{appendix A.1}. With Lemma \ref{lemma 2.1 ODE of Hydrodynamic}, we can state our first main result, which gives the hydrodynamic limit of $\{X_t^N\}_{t\geq 0}$.

\begin{theorem}\label{theorem 2.2 hydrodynamic limit}
Under Assumption (A), for any $t\geq 0$ and $f\in C[0, 1]$,
\[
\lim_{N\rightarrow+\infty}\mu_t^N(f)=\mu_t(f)=\int_0^1\rho(t, u)f(u)du
\]
in $L^2$, where $\{\mu_t\}_{t\geq 0}$ is the unique solution to Equation \eqref{equ definition of hydro ODE} and $\{\rho(t, \cdot)\}_{t\geq 0}$ is the unique solution to Equation \eqref{equ definition of hydro density}.

\end{theorem}

Theorem \ref{theorem 2.2 hydrodynamic limit} shows that the hydrodynamic limit of our linear system $\{X_t^N\}_{t\geq 0}$ is driven by a $\left(C[0, 1]\right)^\prime$-valued ordinary differential equation. An analogue result for the $N$-urn Ehrenfest model is given in Theorem 2.3 of \cite{Xue2022}.

Our second main result is about the fluctuation of $\{X_t^N\}_{t\geq 0}$, i.e., the central limit theorem from the law of large numbers given in Theorem \ref{theorem 2.2 hydrodynamic limit}. We introduce some notations and definitions as a preliminary. For any $t\geq 0$ and $N\geq 1$, we define
\[
V_t^N(du)=\frac{1}{\sqrt{N}}\sum_{i=1}^N\left(X_t^N(i)-\mathbb{E}X_t^N(i)\right)\delta_{\frac{i}{N}}(du),
\]
where $\mathbb{E}$ is the expectation operator. $V_t^N$ is called the fluctuation density field of $\{X_t^N\}_{t\geq 0}$. For given $T>0$, $\{V_t^N\}_{0\leq t\leq T}$ is a random element in $\mathcal{D}\left([0, T], \left(C[0, 1]\right)^\prime\right)$, where $\mathcal{D}\left([0, T], \left(C[0, 1]\right)^\prime\right)$ is the Skorokhod space of c\`{a}dl\`{a}g functions from $[0, T]$ to $\left(C[0, 1]\right)^\prime$ under $\left(C[0, 1]\right)^\prime$ endowed with the weak topology. That is to say, all $F\in \mathcal{D}\left([0, T], \left(C[0, 1]\right)^\prime\right)$ satisfy the following two properties.

(1) For any $0\leq s<T$,
\[
\lim_{t\downarrow s}F_t(f)=F_s(f)
\]
for any $f\in C[0, 1]$.

(2) For any $0<s\leq T$, there exists $F_{s-}\in \left(C[0,1]\right)^\prime$ such that
\[
\lim_{t\uparrow s}F_t(f)=F_{s-}(f)
\]
for any $f\in C[0, 1]$.

For any $t\geq 0$ and $N\geq 1$, we define
\[
\theta_t^N(du)=\frac{1}{N}\sum_{i=1}^N\left(X_t^N(i)\right)^2\delta_{\frac{i}{N}}(du)
\]
and
\[
\omega_t^N(d(u, v))=\frac{1}{N^2}\sum_{i=1}^N\sum_{j=1}^NX_t^N(i)X_t^N(j)\delta_{\left(\frac{i}{N}, \frac{j}{N}\right)}\left(d(u, v)\right).
\]
To give the limit theorem about $\{V_t^N\}_{0\leq t\leq T}$, we first need two laws of large numbers of $\{\theta_t^N\}_{t\geq 0}$ and $\{\omega_t^N\}_{t\geq 0}$. To state these two laws of large numbers, we define $P_3, P_4, P_5, P_6, P_7$ as linear operators from $C[0, 1]$ to $C[0, 1]$ and $l_{1,t}\in \left(C[0, 1]\right)^\prime$ such that
\begin{align*}
&P_3f(u)=\left(c^2(u)-1\right)b(u)f(u),\\
&P_4f(u)=f(u)\int_0^1 \lambda(u, v)\left(a_1^2(u,v)-1\right)dv,\\
&P_5f(u)=\int_0^1\lambda(v, u)a_2^2(v,u)f(v)dv,\\
&P_6f(u)=\int_0^1\lambda(u, v)a_3^2(u,v)f(v)dv,\\
&P_7f(u)=f(u)\int_0^1\lambda(v, u)\left(a_4^2(v, u)-1\right)dv,
\end{align*}
and
\begin{align*}
l_{1, t}(f)=&2\int_0^1\int_0^1\rho(t,u)\rho(t,v)a_1(u,v)a_2(u,v)\lambda(u,v)f(u)dudv\\
&+2\int_0^1\int_0^1\lambda(v, u)a_3(v, u)a_4(v, u)\rho(t, u)\rho(t, v)f(u)dudv
\end{align*}
for any $u\in [0, 1]$ and $f\in C[0, 1]$, where $\rho(t, \cdot)$ is defined as in Lemma \ref{lemma 2.1 ODE of Hydrodynamic}. We further define $P_2=\sum_{k=3}^7P_k$. Now we can state laws of large numbers of $\{\theta_t^N\}_{0\leq t\leq N}$ and $\{\omega_t^N\}_{t\geq 0}$.

\begin{lemma}\label{lemma 2.3 LLNofVartheta}
For any $H\in C\left([0, T]\times[0, 1]\right)$,
\[
\lim_{N\rightarrow+\infty}\sup_{0\leq t\leq T}\mathbb{E}\left(\left(\theta_t^N(H_t)-\theta_t(H_t)\right)^2\right)=0,
\]
where $\{\theta_t\}_{0\leq t\leq T}$ is the unique solution to the $\left(C[0, 1]\right)^\prime$-valued ordinary differential equation
\begin{equation}\label{equ ODE of theta}
\begin{cases}
&\frac{d}{dt}\theta_t=P_2^{*}\theta_t+l_{1,t}\text{~for~}0\leq t\leq T,\\
&\theta_0(du)=\phi(u)du
\end{cases}
\end{equation}
under the norm given in \eqref{equ functionDualInifinityNorm}. Furthermore, there exists $\vartheta\in C\left([0, T]\times[0, 1]\right)$ such that
$\theta_t(du)=\vartheta(t, u)du$ for $0\leq t\leq T$.

\end{lemma}

\begin{lemma}\label{lemma 2.4 LLNofOmega}
For any $t\geq 0$ and $H\in C([0, T]\times[0, 1]\times [0, 1])$,
\[
\lim_{N\rightarrow+\infty}\sup_{0\leq t\leq T}\mathbb{E}\left(\left(\omega_t^N(H_t)-\int_0^1\int_0^1\rho(t,u)\rho(t,v)H_t(u, v)dudv\right)^2\right)=0,
\]
where $\rho(t,\cdot)$ is defined as in Lemma \ref{lemma 2.1 ODE of Hydrodynamic}.
\end{lemma}
Lemmas \ref{lemma 2.3 LLNofVartheta} and \ref{lemma 2.4 LLNofOmega} are called replacement lemmas in investigations of hydrodynamics. Roughly speaking, by these lemmas we can replace the random element $\theta_t^N$ (resp. $\omega_t^N(du, dv)$) by deterministic element $\theta_t$ (resp. $\rho(t,u)\rho(t,v)dudv$) in calculations with errors vanished as $N\rightarrow+\infty$.

We prove Lemmas \ref{lemma 2.3 LLNofVartheta} and \ref{lemma 2.4 LLNofOmega} in Section \ref{section four}. Let $\vartheta$ be defined as in Lemma \ref{lemma 2.3 LLNofVartheta}, then we define
\begin{align*}
K_1(s,u)=&\vartheta(s,u)b(u)\left(c(u)-1\right)^2+\vartheta(s, u)\int_0^1\lambda(u,v)\left(a_1(u,v)-1\right)^2dv\\
&+\int_0^1\vartheta(s,v)\lambda(u,v)a_2^2(u,v)dv+\int_0^1\lambda(v,u)a_3^2(v,u)\vartheta(s,v)dv\\
&+\int_0^1\lambda(v,u)\left(a_4(v,u)-1\right)^2dv\vartheta(s,u)\\
&+2\int_0^1\lambda(u,v)\left(a_1(u,v)-1\right)a_2(u,v)\rho(s,v)dv\rho(s,u)\\
&+2\int_0^1\lambda(v,u)a_3(v,u)\left(a_4(v,u)-1\right)\rho(s,v)dv\rho(s,u)
\end{align*}
and
\begin{align*}
K_2(s,u,v)=&2\lambda(u,v)\left(a_1(u,v)-1\right)a_3(u,v)\vartheta(s,u)\\
&+2\lambda(u,v)a_2(u,v)\left(a_4(u,v)-1\right)\vartheta(s,v)\\
&+2\lambda(u,v)\left(a_1(u,v)-1\right)\left(a_4(u,v)-1\right)\rho(s, u)\rho(s, v)\\
&+2\lambda(u,v)a_2(u,v)a_3(u,v)\rho(s,u)\rho(s,v)
\end{align*}
for any $0\leq s\leq T$ and $u,v\in [0, 1]$. Furthermore, for any $f,g\in C[0, 1]$ and $0\leq s\leq T$, we define
\[
[f,f]_s=\int_0^1K_1(s,u)f^2(u)du+\int_0^1\int_0^1K_2(s,u,v)f(u)f(v)dudv
\]
and
\[
[f,g]_s=\frac{[f+g, f+g]_s-[f-g, f-g]_s}{4}.
\]
Then, we have the following lemma.

\begin{lemma}\label{lemma 2.5}
For any $0\leq s\leq T$, $[\cdot,\cdot]_s$ is a non-negative definite quadratic form, i.e., $[f_1,f_1]_s\geq 0$, $[f_1,f_2]_s=[f_2,f_1]_s$ and $[c_1f_1+c_2f_2,g]_s=c_1[f_1, g]_s+c_2[f_2, g]_s$ for any $f_1, f_2, g\in C[0, 1], c_1, c_2\in \mathbb{R}$.
\end{lemma}
We prove Lemma \ref{lemma 2.5} in Section \ref{section four}. According to Lemma \ref{lemma 2.5}, it is reasonable to define $\{C_t\}_{0\leq t\leq T}$ as the unique stochastic element in $\mathcal{D}\left([0, T], \left(C[0, 1]\right)^\prime\right)$ such that $\{C_t(f)\}_{0\leq t\leq T}$ is a continuous zero-mean martingale with
\[
{\rm Cov}\left(C_t(f), C_t(f)\right)=\int_0^t [f, f]_sds
\]
for any $f\in C[0, 1]$ and $0\leq t\leq T$. According to an analysis similar with that in the proof of Theorem 1.4 of \cite{Holley1978}, there exists a unique stochastic element $V=\{V_t\}_{0\leq t\leq T}$ in $\mathcal{D}\left([0, T], \left(C[0, 1]\right)^\prime\right)$ satisfying all the three following properties.

(1) $\{V_t(f)\}_{0\leq t\leq T}$ is a real-valued continuous function for any $f\in C[0, 1]$.

(2) For any $f\in C[0,1]$ and any $G\in C_c^{\infty}(\mathbb{R})$,
\begin{align*}
\Bigg\{G\left(V_t(f)\right)-G\left(V_0(f)\right)&-\int_0^tG^\prime\left(V_s(f)\right)V_s\left(P_1f\right)ds \\
&-\frac{1}{2}\int_0^tG^{\prime\prime}\left(V_s(f)\right)[f,f]_sds\Bigg\}_{0\leq t\leq T}
\end{align*}
is a martingale.

(3) For any $f\in C[0, 1]$,
\[
V_0(f) \text{~follows~}\mathbb{N}\left(0, \int_0^1f^2(u)\phi(u)(1-\phi(u)) du\right),
\]
where $\mathbb{N}(\mu ,\sigma^2)$ is the normal distribution with mean $\mu$ and variance $\sigma^2$.

Then, it is natural to define the above $V$ as the solution to the $\left(C[0, 1]\right)^\prime$-valued equation
\begin{equation}\label{equ 2.6 generalized O-U}
\begin{cases}
&dV_t=P_1^{*}V_tdt+dC_t\text{~for~}0\leq t\leq T,\\
&V_0(f) \text{~follows~}\mathbb{N}\left(0, \int_0^1f^2(u)\phi(u)(1-\phi(u))du\right)\text{~for any~}f\in C[0, 1],\\
&V_0 \text{~is independent of~}\{C_t\}_{0\leq t\leq T}.
\end{cases}
\end{equation}
Now we can state our second main result, which gives the fluctuation of $\{X_t^N\}_{t\geq 0}$.

\begin{theorem}\label{theorem 2.6 fluctuation of the NurnLinearSyetem}
Let $V^N=\{V_t^N\}_{0\leq t\leq T}$ and $V=\{V_t\}_{0\leq t\leq T}$ be the solution to \eqref{equ 2.6 generalized O-U}, then $V^N$ converges weakly to $V$ as $N\rightarrow+\infty$.
\end{theorem}

Theorem \ref{theorem 2.6 fluctuation of the NurnLinearSyetem} shows that the fluctuation of our linear system $\{X_t^N\}_{t\geq 0}$ is driven by a $\left(C[0, 1]\right)^\prime$-valued Ornstein-Uhlenbeck process. An analogue result for the $N$-urn Ehrenfest model is given in Theorem 2.4 of \cite{Xue2022}.

For any $t\in \mathbb{R}$, it is reasonable to define
\[
e^{tP_1}=\sum_{k=0}^{+\infty}\frac{t^kP_1^k}{k!}\text{~and~}e^{tP_1^{*}}=\sum_{k=0}^{+\infty}\frac{t^k \left(P_1^{*}\right)^k}{k!}
\]
since
\[
\|P_1f\|_\infty\leq \left(\|b\|_\infty\left(\|c\|_\infty+1\right)+\|\lambda\|_\infty\left(2+\sum_{i=1}^4\|a_i\|_\infty\right)\right)\|f\|_\infty
\]
and
\[
\|P_1^{*}\nu\|\leq \left(\|b\|_\infty\left(\|c\|_\infty+1\right)+\|\lambda\|_\infty\left(2+\sum_{i=1}^4\|a_i\|_\infty\right)\right)\|\nu\|
\]
for any $f\in C[0, 1]$ and $\nu\in \left(C[0, 1]\right)^\prime$.
According to an analysis similar with those in the proof of Theorem 1.4 of \cite{Holley1978} and Appendix A.2 of \cite{Xue2022}, for any $t>s$ and $f\in C[0, 1]$,
\[
\mathbb{E}\left(e^{\sqrt{-1}V_t(f)}\big|V_r, r\leq s\right)=e^{\sqrt{-1}V_s\left(e^{(t-s)P_1}f\right)-\frac{1}{2}\int_s^t\left[e^{(t-r)P_1}f,e^{(t-r)P_1}f\right]_rdr}.
\]
As a result, for any $t\geq 0$ and $f\in C[0, 1]$, $V_t(f)$ follows
\begin{equation}\label{equ 2.7}
\mathbb{N}\left(0, \int_0^1\left(\left(e^{tP_1}f\right)(u)\right)^2\phi(u)(1-\phi(u))du+\int_0^t \left[e^{(t-s)P_1}f, e^{(t-s)P_1}f\right]_sds\right).
\end{equation}
Another way to deduce \eqref{equ 2.7} is solving Equation \eqref{equ 2.6 generalized O-U} directly. In detail, Equation \eqref{equ 2.6 generalized O-U} implies that $d(e^{-tP_1^*}V_t)=e^{-tP_1^*}dC_t$ and hence
\[
V_t=e^{tP_1^*}V_0+\int_0^te^{(t-s)P_1^*}dC_s,
\]
which also leads to \eqref{equ 2.7}. In conclusion, we have the following corollary of Theorem \ref{theorem 2.6 fluctuation of the NurnLinearSyetem}.

\begin{corollary}\label{corollary 2.7}
For any $t\geq 0$ and $f\in C[0, 1]$, $V_t^N(f)$ converges in distribution to
\[
\mathbb{N}\left(0, \int_0^1\left(\left(e^{tP_1}f\right)(u)\right)^2\phi(u)(1-\phi(u))du+\int_0^t \left[e^{(t-s)P_1}f, e^{(t-s)P_1}f\right]_sds\right)
\]
as $N\rightarrow+\infty$.
\end{corollary}

The proof of Theorem \ref{theorem 2.2 hydrodynamic limit} is given in Section \ref{section three}. We first show that
\[
\lim_{N\rightarrow+\infty}\mathbb{E}\mu_t^N(f)=\mu_t(f)
\]
and then show that
\[
\lim_{N\rightarrow+\infty}{\rm Var}\left(\mu_t^N(f)\right)=0.
\]
The proof of Theorem \ref{theorem 2.6 fluctuation of the NurnLinearSyetem} is given in Section \ref{section four}. We first prove replacement lemmas \ref{lemma 2.3 LLNofVartheta} and \ref{lemma 2.4 LLNofOmega}, then by Dynkin's martingale formula and Lemmas \ref{lemma 2.3 LLNofVartheta}, \ref{lemma 2.4 LLNofOmega}, we show that any weak limit of a subsequence of $\{V^N\}_{N\geq 1}$ is the solution to the martingale problem related to the generalized O-U process in \eqref{equ 2.6 generalized O-U}. An extension in linear systems of Chapman-Kolmogorov equation, which follows from theorems given in Chapter 9 of \cite{Lig1985}, plays key role in the above strategy. 

\section{Proof of Theorem \ref{theorem 2.2 hydrodynamic limit}}\label{section three}
In this section we prove Theorem \ref{theorem 2.2 hydrodynamic limit}. For later use, we first introduce some notations. For any $N\geq 1$, we use $\{S^N(t)\}_{t\geq 0}$ to denote the Markov semigroup of $\{X_t^N\}_{t\geq 0}$, i.e.,
\[
S^N(t)H(x)=\mathbb{E}\left(H(X_t^N)\big|X_0^N=x\right)
\]
for any $x\in [0, +\infty)^N$ and $H$ from $[0, +\infty)^N$ to $\mathbb{R}$. For any $t\geq 0$ and $N\geq 1$, we use $\hat{\mu}_t^N$ to denote the element in $\left(C[0, 1]\right)^\prime$ such that
\[
\hat{\mu}_t^N(f)=\mathbb{E}\mu_t^N(f)
\]
for any $f\in C[0, 1]$ and use $\hat{\omega}_t^N, \varpi_t^N$ to denote elements in $\left(C\left([0, 1]\times[0, 1]\right)\right)^\prime$ such that
\[
\hat{\omega}_t^N(H)=\mathbb{E}\omega_t^N(H)=\frac{1}{N^2}\sum_{i=1}^N\sum_{j=1}^N\mathbb{E}\left(X_t^N(i)X_t^N(j)\right)H\left(\frac{i}{N}, \frac{j}{N}\right)
\]
and
\[
\varpi_t^N(H)=\frac{1}{N^2}\sum_{i=1}^N\sum_{j=1}^N\mathbb{E}\left(X_t^N(i)\right)\mathbb{E}\left(X_t^N(j)\right)H\left(\frac{i}{N}, \frac{j}{N}\right)
\]
for any $H\in C\left([0, 1]\times [0, 1]\right)$. Note that $\hat{\mu}_t^N, \hat{\omega}_t^N, \varpi_t^N$ are all deterministic.

Theorem \ref{theorem 2.2 hydrodynamic limit} follows from the following two lemmas.

\begin{lemma}\label{lemma 3.1}
For any $t\geq 0$ and $f\in C[0, 1]$,
\[
\lim_{N\rightarrow+\infty}\hat{\mu}_t^N(f)=\mu_t(f),
\]
where $\mu=\{\mu_t\}_{t\geq 0}$ is the solution to Equation \eqref{equ definition of hydro ODE}.
\end{lemma}

\begin{lemma}\label{lemma 3.2}
For any $t\geq 0$ and $f\in C[0, 1]$,
\[
\lim_{N\rightarrow+\infty}{\rm Var}\left(\mu_t^N(f)\right)=0.
\]
\end{lemma}

\proof[Proof of Theorem \ref{theorem 2.2 hydrodynamic limit}]

Since $(c+d)^2\leq 2c^2+2d^2$ for any $c,d\in \mathbb{R}$,
\[
\mathbb{E}\left(\left(\mu_t^N(f)-\mu_t(f)\right)^2\right)\leq 2{\rm Var}\left(\mu_t^N(f)\right)+2\left(\hat{\mu}_t^N(f)-\mu_t(f)\right)^2
\]
and hence Theorem \ref{theorem 2.2 hydrodynamic limit} holds according to Lemmas \ref{lemma 3.1} and \ref{lemma 3.2}.

\qed

The remainder of this section is devoted to proofs of Lemmas \ref{lemma 3.1} and \ref{lemma 3.2}, where we will frequently execute the calculation
\begin{equation}\label{equ linear system Chapman-Kolmogorov}
\frac{d}{dt}S^N(t)H=S^N(t)\mathcal{L}_NH
\end{equation}
for $H$ from $[0, +\infty)^N$ to $\mathbb{R}$ with the form $H(x)=\prod_{j=1}^kx(l_j)$ for any $x\in [0, +\infty)^N$ and some integers $k\geq 1, 1\leq l_1, l_2,\ldots, l_k\leq N$. Such an extension of Chapman-Kolmogorov equation holds according to Theorems 9.1.27 and 9.3.1 of \cite{Lig1985}.

Now we prove Lemma \ref{lemma 3.1}

\proof

By \eqref{equ linear system Chapman-Kolmogorov} and the definition of $\mathcal{L}_N$ given Section \ref{section one},
\begin{align}\label{equ 3.2}
&\frac{d}{dt}\mathbb{E}X_t^N(i)= \notag\\
&b\left(\frac{i}{N}\right)\left(c\left(\frac{i}{N}\right)-1\right)\mathbb{E}X_t^N(i)
+\frac{\mathbb{E}X_t^N(i)}{N}\sum_{j\neq i}\lambda\left(\frac{j}{N}, \frac{i}{N}\right)\left(a_4\left(\frac{j}{N}, \frac{i}{N}\right)-1\right)\notag\\
&+\frac{\mathbb{E}X_t^N(i)}{N}\sum_{j\neq i}\lambda\left(\frac{i}{N}, \frac{j}{N}\right)\left(a_1\left(\frac{i}{N}, \frac{j}{N}\right)-1\right)\notag\\
&+\frac{1}{N}\sum_{j\neq i}\lambda\left(\frac{i}{N}, \frac{j}{N}\right)a_2\left(\frac{i}{N}, \frac{j}{N}\right)\mathbb{E}X_t^N(j)\\
&+\frac{1}{N}\sum_{j\neq i}\lambda\left(\frac{j}{N}, \frac{i}{N}\right)a_3\left(\frac{j}{N}, \frac{i}{N}\right)\mathbb{E}X_t^N(j) \notag
\end{align}
for each $1\leq i\leq N$. Therefore, there exists $C_2=C_2(b, c, \lambda, a_1, a_2, a_3, a_4)<+\infty$ independent of $N, t$ such that
\[
\frac{d}{dt}\left(\frac{1}{N}\sum_{i=1}^N\mathbb{E}X_t^N(i)\right)\leq C_2 \left(\frac{1}{N}\sum_{i=1}^N\mathbb{E}X_t^N(i)\right)
\]
and hence $\frac{1}{N}\sum_{i=1}^N\mathbb{E}X_t^N(i)\leq \|\phi\|_\infty e^{C_2 t}$ for any $t\geq 0$ according to Assumption (A). Then, for any $T_1>0$, there exists $C_3=C_3(T_1)<+\infty$ independent of $t, i, N$ such that
\[
\frac{d}{dt}\mathbb{E}X_t^N(i)\leq C_3\mathbb{E}X_t^N(i)+C_3
\]
for all $0\leq t\leq T_1, N\geq 1$ and $1\leq i\leq N$. As a result, there exists $C_4=C_4(T_1)<+\infty$ independent of $t, i, N$ such that
\[
\frac{d}{dt}\mathbb{E}X_t^N(i), \mathbb{E}X_t^N(i)\leq C_4
\]
for all $0\leq t\leq T_1, N\geq 1$ and $1\leq i\leq N$. Hence, for any given $f\in C[0, 1]$, $\{\hat{\mu}_t^N(f):~0\leq t\leq T_1\}_{N\geq 1}$ are uniformly bounded and equicontinuous. Then, by Ascoli-Arzela theorem, any subsequence $\{\hat{\mu}_t^{N_k}:~0\leq t\leq T_1\}_{k\geq 1}$ of $\{\hat{\mu}_t^N:~0\leq t\leq T_1\}_{N\geq 1}$ has a subsequence $\{\hat{\mu}_t^{N_{k_j}}:~0\leq t\leq T_1\}_{j\geq 1}$ such that
\[
\lim_{j\rightarrow+\infty}\sup_{0\leq t\leq T_1}\left|\hat{\mu}_t^{N_{k_j}}(f)-\nu_t(f)\right|=0
\]
for some $\nu\in \mathcal{D}\left([0, T_1], \left(C[0, 1]\right)^\prime\right)$ and any $f\in C[0, 1], 0\leq t\leq T_1$. To complete this proof, we only need to show that $\nu$ equals $\{\mu_t\}_{0\leq t\leq T_1}$. By Assumption (A), Equation \eqref{equ 3.2} and the fact that $\sup_{0\leq t\leq T_1}\mathbb{E}X_t^N(i)\leq C_4$,
\begin{equation}\label{equ 3.3}
\hat{\mu}_t^{N_{k_j}}(f)=\frac{1}{N_{k_j}}\sum_{i=1}^{N_{k_j}}\phi\left(\frac{i}{N_{k_j}}\right)f\left(\frac{i}{N_{k_j}}\right)
+\int_0^t \hat{\mu}_s^{N_{k_j}}\left(P_{1,{N_{k_j}}}f\right)ds+O\left(\frac{1}{N_{k_j}}\right)
\end{equation}
for any $f\in C[0, 1]$ and $0\leq t\leq T_1$, where
\begin{align*}
&P_{1,N}f(u)=\\
&\frac{f(u)}{N}\sum_{j=1}^N\lambda\left(u, \frac{j}{N}\right)\left(a_1\left(u, \frac{j}{N}\right)-1\right)+b(u)\left(c(u)-1\right)f(u)\\
&+\frac{1}{N}\sum_{j=1}^N\lambda\left(u, \frac{j}{N}\right)a_3\left(u, \frac{j}{N}\right)f\left(\frac{j}{N}\right)\\
&+\frac{f(u)}{N}\sum_{j=1}^N\lambda\left(\frac{j}{N}, u\right)\left(a_4\left(\frac{j}{N}, u\right)-1\right)\\
&+\frac{1}{N}\sum_{j=1}^N\lambda\left(\frac{j}{N}, u\right)f\left(\frac{j}{N}\right)a_2\left(\frac{j}{N}, u\right)
\end{align*}
for any $u\in [0, 1]$. According to the definition of $P_{1, N}$,
\begin{equation}\label{equ 3.4}
\lim_{N\rightarrow+\infty}\|P_{1, N}f-P_1f\|_\infty=0
\end{equation}
for any $f\in C[0, 1]$. According to the fact that $\sup_{0\leq t\leq T_1}\mathbb{E}X_t^N(i)\leq C_4$,
\begin{equation}\label{equ 3.5}
\sup_{0\leq t\leq T_1}\|\mu_t^N\|\leq C_4.
\end{equation}
By \eqref{equ 3.3}, \eqref{equ 3.4} and \eqref{equ 3.5}, let $j\rightarrow+\infty$, then
\begin{equation}\label{equ 3.6}
\nu_t(f)=\int_0^1\phi(u)f(u)du+\int_0^t\nu_s(P_1f)ds
\end{equation}
for all $f\in C[0, 1]$ and $0\leq t\leq T_1$. As we have shown in the proof of Lemma \ref{lemma 2.1 ODE of Hydrodynamic}, Equation \eqref{equ 3.6} implies that
$\nu=\{\mu_t\}_{0\leq t\leq T_1}$ according to Gronwall's inequality and hence the proof is complete.

\qed

At last we only need to prove Lemma \ref{lemma 3.2}. As a preliminary, we introduce some notations and definitions. For each $N\geq 1$, we use $G_N$ to denote  the set $\{1,2,\ldots,N\}\times \{1,2,\ldots, N\}$. For any $t\geq 0$, we define $F_t^N$ and $\hat{F}_t^N$ as functions from $G_N$ to $[0, +\infty)$ such that
\[
F_t^N(i, j)=\mathbb{E}\left(X_t^N(i)\right)\mathbb{E}\left(X_t^N(j)\right)
\]
and
\[
\hat{F}_t^N(i,j)=\mathbb{E}\left(X_t^N(i)X_t^N(j)\right)
\]
for any $1\leq i, j\leq N$. For a finite set $S$, a function from $S\times S$ to $\mathbb{R}$ is called a $S\times S$ matrix. For any $S\times S$ matrices $M_1, M_2$, function $H$ from $S$ to $\mathbb{R}$ and $t\geq 0$, $M_1M_2$ is defined as the $S\times S$ matrix such that
\[
M_1M_2(x,y)=\sum_{z\in S}M_1(x,z)M_2(z,y)
\]
for any $x,y\in S$ and $M_1H$ is defined as the function from $S$ to $\mathbb{R}$ such that
\[
M_1H(x)=\sum_{y\in S}M_1(x,y)H(y)
\]
and $e^{tM_1}$ is defined as the $S\times S$ matrix such that
\[
e^{tM_1}=\sum_{k=0}^{+\infty}\frac{t^kM_1^k}{k!}.
\]
Now we prove Lemma \ref{equ 3.2}.

\proof{Proof of Lemma \ref{lemma 3.2}}

By \eqref{equ linear system Chapman-Kolmogorov}, there exist $G_N\times G_N$ matrices $\mathcal{M}_N$ and $\hat{\mathcal{M}}_N$ such that
\[
\frac{d}{dt}F_t^N=\mathcal{M}_NF_t^N \text{~and~} \frac{d}{dt}\hat{F}_t^N=\hat{\mathcal{M}}_N\hat{F}_t^N
\]
for any $t\geq 0$. The expressions of $\mathcal{M}_N$ and $\hat{\mathcal{M}}_N$ are a little tedious, which we put in Appendix \ref{appendix A.2}. Then,
\begin{equation}\label{equ 3.7}
F_t^N=e^{t\mathcal{M}_N}F_0^N \text{~and~}\hat{F}_t^N=e^{t\hat{\mathcal{M}}_N}\hat{F}_0^N.
\end{equation}
According to definitions of $\mathcal{M}_N$ and $\hat{\mathcal{M}}_N$ given in Appendix \ref{appendix A.2}, there exist $C_5=C_5\left(b, c, \lambda, a_1, a_2, a_3, a_4\right)<+\infty$ independent of $N$ such that
\[
\sup_{(i, j)\in G_N}\left|\mathcal{M}_NH(i, j)\right|, \sup_{(i, j)\in G_N}\left|\hat{\mathcal{M}}_NH(i, j)\right|\leq C_5\sup_{(i, j)\in G_N}\left|H(i,j)\right|
\]
for any function $H$ from $G_N$ to $\mathbb{R}$. As a result, by Assumption (A),
\begin{equation}\label{equ 3.8}
F_t^N(i,j), \hat{F}_t^N(i,j)\leq e^{C_5t}\|\phi\|_\infty
\end{equation}
for any $t\geq 0, N\geq 1, 1\leq i,j\leq N$. For any $t\geq 0$, we claim that
\begin{equation}\label{equ 3.9}
\lim_{N\rightarrow+\infty}\sup_{1\leq i\neq j\leq N}\left|\hat{F}_t^N(i,j)-F_t^N(i,j)\right|=0.
\end{equation}
We check \eqref{equ 3.9} at the end of this proof. By \eqref{equ 3.8},
\begin{align*}
{\rm Var}\left(\mu_t^N(f)\right)&=\frac{1}{N^2}\sum_{i=1}^N\sum_{j=1}^N\left(\hat{F}_t^N(i,j)-F_t^N(i,j)\right)f\left(\frac{i}{N}\right)f\left(\frac{j}{N}\right)\\
&\leq \frac{2}{N}e^{C_5t}\|f\|_\infty^2+\left(\sup_{1\leq i\neq j\leq N}\left|\hat{F}_t^N(i,j)-F_t^N(i,j)\right|\right)\|f\|_\infty^2
\end{align*}
and hence Lemma \ref{lemma 3.2} holds according to \eqref{equ 3.9}.

At last we check \eqref{equ 3.9}. For integer $r\geq 0$, let
\[
\alpha_r^N=\sup_{1\leq i\neq j\leq N}\left|\mathcal{M}_N^rF_0^N(i,j)-\hat{\mathcal{M}}_N^r\hat{F}_0^N(i,j)\right|,
\]
then $a_0^N=0$ according to Assumption (A). For all $r\geq 0$ and $i\neq j$,
\begin{align}\label{equ 3.10}
&\hat{\mathcal{M}}_N^{r+1}\hat{F}_0^N(i,j)-M_N^{r+1}F_0^N(i,j)\\
&=\sum_{l,k}\left(\hat{\mathcal{M}}_N\left((i,j), (l,k)\right)\hat{\mathcal{M}}_N^r\hat{F}_0^N(l,k)
-\mathcal{M}_N\left((i,j), (l,k)\right)\mathcal{M}_N^rF_0^N(l,k)\right). \notag
\end{align}
When $i\neq j$, $l\neq k$ and $(l,k)\neq (i,j), (i,i), (j,j)$, by definitions of $\mathcal{M}_N$ and $\hat{\mathcal{M}}_N$, there exists $C_6=C_6(b, c, \lambda, a_1, a_2, a_3, a_4)<+\infty$ independent of $N, i, j, l, k$ such that
\begin{equation}\label{equ 3.13}
0\leq \hat{\mathcal{M}}_N\left((i,j), (l,k)\right)=M_N\left((i,j), (l,k)\right)\leq \frac{C_6}{N}.
\end{equation}
When $i\neq j$ and $(l,k)=(i,j)$, there exist $C_8= C_8(b, c, \lambda, a_1, a_2, a_3, a_4)<+\infty$ and $C_9=C_9(b, c, \lambda, a_1, a_2, a_3, a_4)<+\infty$ independent of $N, i, j$ such that
\begin{equation}\label{equ 3.11}
|\mathcal{M}_N\left((i,j), (i,j)\right)|\leq C_8 \text{~and~} |\hat{\mathcal{M}}_N\left((i,j), (i,j)\right)-\mathcal{M}_N\left((i,j), (i,j)\right)|\leq \frac{C_9}{N}.
\end{equation}
When $i\neq j$ and $(l,k)=(i, i)$ or $(j, j)$, there exist $C_{10}=C_{10}(b, c, \lambda, a_1, a_2, a_3, a_4)<+\infty$ independent of $N, i, j$ such that
\begin{equation}\label{equ 3.12}
|\mathcal{M}_N\left((i,j), (l, k)\right)|, |\hat{\mathcal{M}}_N\left((i,j), (l, k)\right)|\leq \frac{C_{10}}{N}.
\end{equation}
By \eqref{equ 3.10}, \eqref{equ 3.13}, \eqref{equ 3.11} and \eqref{equ 3.12},
\begin{align*}
a_{r+1}^N&\leq \frac{2C_6N}{N}a_r^N+\frac{C_9}{N}C_5^r\|\phi\|_\infty+C_8a_r^N+\frac{4C_{10}}{N}C_5^r\|\phi\|_\infty\\
&=C_7a_r^N+\frac{C_{11}C_5^r\|\phi\|_\infty}{N}
\end{align*}
for all $r\geq 0$, where $C_7=2C_6+C_8$ and $C_{11}=4C_{10}+C_9$. Since $a_0^N=0$, by induction,
\begin{equation}\label{equ 3.14}
a_r^N\leq \frac{C_{11}\|\phi\|_\infty}{N}\left(\sum_{l=0}^{r-1}C_7^{r-1-l}C_5^l\right)\leq \frac{C_{11}\|\phi\|_\infty}{N}(C_7+C_5)^{r-1}
\end{equation}
for all $r\geq 1$. By \eqref{equ 3.7}, for $i\neq j$,
\begin{equation}\label{equ 3.15}
\left|\hat{F}_t^N(i,j)-F_t^N(i,j)\right|\leq \sum_{r=0}^{+\infty}\frac{t^ra_r^N}{r!}\leq \frac{C_{11}\|\phi\|_\infty te^{(C_7+C_5)t}}{N}
\end{equation}
and hence \eqref{equ 3.9} follows from \eqref{equ 3.15}.

\qed

\section{Proof of Theorem \ref{theorem 2.6 fluctuation of the NurnLinearSyetem}}\label{section four}
In this section we prove Theorem \ref{theorem 2.6 fluctuation of the NurnLinearSyetem}. As preliminaries, we first prove Lemmas \ref{lemma 2.4 LLNofOmega} and \ref{lemma 2.3 LLNofVartheta} in Subsections \ref{subsection 4.1} and \ref{subsection 4.2} respectively. According to these two lemmas and Dynkin's martingale formula, we prove Theorem \ref{theorem 2.6 fluctuation of the NurnLinearSyetem} in Subsection \ref{subsection 4.3} by showing that any weak limit of $\{V^N\}_{N\geq 1}$ is the solution to the martingale problem related to the generalized O-U process given in \eqref{equ 2.6 generalized O-U}.

\subsection{Proof of Lemma \ref{lemma 2.4 LLNofOmega}}\label{subsection 4.1}
In this subsection we give the proof of Lemma \ref{lemma 2.4 LLNofOmega}, which follows an analysis similar with that in the proof of Theorem \ref{theorem 2.2 hydrodynamic limit}.

\proof[Proof of Lemma \ref{lemma 2.4 LLNofOmega}]

Let $\hat{\omega}_t^N$ be defined as in Section \ref{section three}, then we claim that
\begin{equation}\label{equ 4.1}
\lim_{N\rightarrow+\infty}\sup_{0\leq t\leq T}\left|\hat{\omega}_t^N(H_t)-\int_0^1\int_0^1\rho(t,u)\rho(t,v)H_t(u,v)dudv\right|=0
\end{equation}
and
\begin{equation}\label{equ 4.2}
\lim_{N\rightarrow+\infty}\sup_{0\leq t\leq T}{\rm Var}\left(\omega_t^N(H_t)\right)=0
\end{equation}
for any $H\in C\left([0, T]\times[0, 1]\times [0, 1]\right)$. We check \eqref{equ 4.1} and \eqref{equ 4.2} later. Since
\begin{align*}
&\mathbb{E}\left(\left(\omega_t^N(H_t)-\int_0^1\int_0^1\rho(t,u)\rho(t,v)H_t(u,v)dudv\right)^2\right)\\
&\leq 2\left(\hat{\omega}_t^N(H_t)-\int_0^1\int_0^1\rho(t,u)\rho(t,v)H_t(u,v)dudv\right)^2+2{\rm Var}\left(\omega_t^N(H_t)\right),
\end{align*}
Lemma \ref{lemma 2.4 LLNofOmega} follows from \eqref{equ 4.1} and \eqref{equ 4.2}.

\qed

The remainder of this subsection is devoted to proofs of \eqref{equ 4.1} and \eqref{equ 4.2}.

\proof[Proof of \eqref{equ 4.1}]

For integer $k, l, r\geq 0$, let $H_{k,l,r}$ be the function from $[0, T]\times[0, 1]\times [0, 1]$ to $\mathbb{R}$ such that $H_{k,l, r}(t, u, v)=t^ru^kv^l$ for all $u, v\in [0, 1], t\in [0, T]$. Since ${\rm span}\{H_{k,l,r}:~k,l,r\geq 0\}$ are dense in $C\left([0, T]\times[0, 1]\times[0, 1]\right)$ and
\[
\sup_{N\geq 1}\sup_{1\leq i\leq N}\sup_{0\leq t\leq T}\mathbb{E}X_t^N(i)<+\infty
\]
as we have shown in Section \ref{section three}, we only need to show that \eqref{equ 4.1} holds for all $H_{k,l,r}$. Let $\hat{\mu}_t^N, \varpi_t^N$ be defined as in Section \ref{section three}, then $\varpi_t^N\left(H_{k,l,r}(t,\cdot)\right)=t^r\hat{\mu}_t^N(\beta_k)\hat{\mu}_t^N(\beta_l)$, where
$\beta_k(u)=u^k$ for all integer $k\geq 0$. Since $\beta_k\in C[0, 1]$, as we have shown in the proof of Lemma \ref{lemma 3.1},
\[
\lim_{N\rightarrow+\infty}\varpi_t^N\left(H_{k,l,r}(t,\cdot)\right)=t^r\mu_t(\beta_k)\mu_t(\beta_l)=\int_0^1\int_0^1\rho(t,u)\rho(t,v)H_{k,l,r}(t,u,v)dudv
\]
uniformly in $[0, T]$. By \eqref{equ 3.8} and \eqref{equ 3.15},
\[
\sup_{0\leq t\leq T}\left|\hat{\omega}_t^N\left(H_{k,l,r}(t,\cdot)\right)-\varpi_t^N\left(H_{k,l,r}(t,\cdot)\right)\right|\leq \frac{2e^{C_5T}\|\phi\|_\infty+C_{11}\|\phi\|_\infty Te^{(C_7+C_5)T}}{N}T^r
\]
and hence \eqref{equ 4.1} holds for all $H_{k,l,r}$.

\qed

\proof[Proof of \eqref{equ 4.2}]

For any $t\geq 0$ and $N\geq 1$, let $\Lambda_t^N$ and $\hat{\Lambda}_t^N$ be functions from $G_N\times G_N$ to $\mathbb{R}$ such that
\[
\Lambda_t^N(i,j, l,k)=\hat{F}_t^N(i,j)\hat{F}_t^N(l,k)=\mathbb{E}\left(X_t^N(i)X_t^N(j)\right)\mathbb{E}\left(X_t^N(l)X_t^N(k)\right)
\]
and
\[
\hat{\Lambda}_t^N(i,j,l,k)=\mathbb{E}\left(X_t^N(i)X_t^N(j)X_t^N(l)X_t^N(k)\right)
\]
for any $1\leq i, j, l, k\leq N$. By \eqref{equ linear system Chapman-Kolmogorov}, there exist $\left(G_N\times G_N\right)\times \left(G_N\times G_N\right)$ matrices $\Xi_N$ and $\hat{\Xi}_N$ such that
\[
\frac{d}{dt}\Lambda_t^N=\Xi_N\Lambda_t^N \text{~and~}\frac{d}{dt}\hat{\Lambda}_t^N=\hat{\Xi}_N\hat{\Lambda}_t^N
\]
and hence $\Lambda_t^N=e^{t\Xi_N}\Lambda_0^N, \hat{\Lambda}_t^N=e^{t\hat{\Xi}_N}\hat{\Lambda}_0^N$ for all $t\geq 0$. Since detailed expressions of $\Xi_N$ and $\hat{\Xi}_N$ are too tedious, we leave them for readers. Here we give some properties of $\Xi_N$ and $\hat{\Xi}_N$ for later use, which are not difficult (but very tedious) to check by \eqref{equ linear system Chapman-Kolmogorov}. There exists $C_{12}=C_{12}\left(b, c, \lambda, a_1, a_2, a_3, a_4\right)<+\infty$ independent of $N$ such that
\begin{align}\label{equ 4.3}
&\sup_{(i,j,l,k)\in G_N\times G_N}\left|\Xi_NH(i,j,l,k)\right|, \sup_{(i,j,l,k)\in G_N\times G_N}\left|\hat{\Xi}_NH(i,j,l,k)\right| \notag \\
&\leq C_{12}\sup_{(i,j,l,k)\in G_N\times G_N}\left|H(i,j,l,k)\right|
\end{align}
for any function $H$ from $G_N\times G_N$ to $\mathbb{R}$. There exists $C_{13}<+\infty$ independent of $N$ such that for all $(i,j,l,k)\in G_N\times G_N$,
\begin{equation}\label{equ 4.4}
{\rm card}\left\{(i_2, j_2, l_2, k_2)\in G_N\times G_N:~\Xi_N\left((i,j,l,k),(i_2,j_2,l_2,k_2)\right)\neq 0\right\}\leq C_{13}N
\end{equation}
and
\begin{equation}\label{equ 4.5}
{\rm card}\left\{(i_2, j_2, l_2, k_2)\in G_N\times G_N:~\hat{\Xi}_N\left((i,j,l,k),(i_2,j_2,l_2,k_2)\right)\neq 0\right\}\leq C_{13}N
\end{equation}
where ${\rm card}(\cdot)$ is the cardinality function. There exists
\[
C_{14}=C_{14}(b, c, \lambda, a_1, a_2, a_3, a_4)<+\infty
\]
independent of $N$ such that for all $(i,j,l,k), (i_2, j_2, l_2, k_2)\in G_N\times G_N$ satisfying $(i,j,l,k)\neq (i_2, j_2, l_2, k_2)$,
\begin{equation}\label{equ 4.6}
|\Xi_N\left((i,j,l,k),(i_2,j_2,l_2,k_2)\right)|, |\hat{\Xi}_N\left((i,j,l,k),(i_2,j_2,l_2,k_2)\right)|\leq \frac{C_{14}}{N}.
\end{equation}
There exists $C_{15}<+\infty$ independent of $N$ such that for any $(i,j,l,k)\in G_N\times G_N$ satisfying $\{i, j\}\cap\{l,k\}=\emptyset$,
\begin{align}\label{equ 4.7}
&{\rm card}\Bigg\{(i_2, j_2, l_2, k_2)\in G_N\times G_N:\notag\\
&\text{\quad\quad}\Xi_N\left((i,j,l,k),(i_2,j_2,l_2,k_2)\right)\neq \hat{\Xi}_N\left((i,j,l,k),(i_2,j_2,l_2,k_2)\right)\Bigg\}\leq C_{15}.
\end{align}
There exist $C_{16}=C_{16}(b, c, \lambda, a_1, a_2, a_3, a_4)<+\infty$ and
\[
C_{17}=C_{17}(b, c, \lambda, a_1, a_2, a_3, a_4)<+\infty
\]
independent of $N$ such that for all $(i,j, l,k)\in G_N\times G_N$ satisfying $\{i, j\}\cap\{l,k\}=\emptyset$,
\begin{equation}\label{equ 4.8}
|\Xi_N\left((i,j,l,k), (i,j,l,k)\right)|\leq C_{16}
\end{equation}
and
\begin{equation}\label{equ 4.9}
\left|\Xi_N\left((i,j,l,k), (i,j,l,k)\right)-\hat{\Xi}_N\left((i,j,l,k), (i,j,l,k)\right)\right|\leq \frac{C_{17}}{N}.
\end{equation}
There exists $C_{18}<+\infty$ independent of $N$ such that for any $(i,j,l,k)\in G_N\times G_N$ satisfying $\{i,j\}\cap\{l,k\}=\emptyset$,
\begin{align}\label{equ 4.9 two}
&{\rm card}\Big\{(i_2,j_2,l_2,k_2)\in G_N\times G_N:~\{i_2, j_2\}\cap\{l_2, k_2\}\neq \emptyset \notag\\
&\text{\quad\quad and~}\Xi_N\left((i,j,l,k),(i_2,j_2, l_2, k_2)\right)\neq 0\Big\}\leq C_{18}
\end{align}
and
\begin{align}\label{equ 4.9 three}
&{\rm card}\Big\{(i_2,j_2,l_2,k_2)\in G_N\times G_N:~\{i_2, j_2\}\cap\{l_2, k_2\}\neq \emptyset \notag\\
&\text{\quad\quad and~}\hat{\Xi}_N\left((i,j,l,k),(i_2,j_2, l_2, k_2)\right)\neq 0\Big\}\leq C_{18}.
\end{align}
By \eqref{equ 4.3} and Assumption (A), for any $t\geq 0, N\geq 1$ and $(i,j,l,k)\in G_N\times G_N$,
\begin{equation}\label{equ 4.10}
\Lambda_t^N(i,j,l,k), \hat{\Lambda}_t^N(i,j,l,k)\leq e^{tC_{12}}.
\end{equation}
We claim that
\begin{equation}\label{equ 4.11}
\lim_{N\rightarrow+\infty}\sup_{0\leq t\leq T}\sup_{(i,j,l,k)\in G_N\times G_N:
\atop \{i,j\}\cap\{l,k\}=\emptyset}\left|\Lambda_t^N(i,j,l,k)-\hat{\Lambda}_t^N(i,j,l,k)\right|=0.
\end{equation}
We check \eqref{equ 4.11} later. By \eqref{equ 4.10}, for $0\leq t\leq T$,
\begin{align*}
&{\rm Var}\left(\omega_t^N(H_t)\right)\\
&=\frac{1}{N^4}\sum_{(i,j,l,k)\in G_N\times G_N}\left(\Lambda_t^N(i,j,l,k)-\hat{\Lambda}_t^N(i,j,l,k)\right)H_t\left(\frac{i}{N}, \frac{j}{N}\right)H_t\left(\frac{l}{N}, \frac{k}{N}\right)\\
&\leq \frac{4N^3}{N^4}e^{tC_{12}}\|H\|^2_{\infty}+\sup_{(i,j,l,k)\in G_N\times G_N:
\atop \{i,j\}\cap\{l,k\}=\emptyset}\left|\Lambda_t^N(i,j,l,k)-\hat{\Lambda}_t^N(i,j,l,k)\right|\|H\|^2_{\infty}
\end{align*}
and hence \eqref{equ 4.2} follows from \eqref{equ 4.11}.

At last we check \eqref{equ 4.11}. For integers $r\geq 0$ and $N\geq 1$, let
\[
h_r^N=\sup_{(i,j,l,k)\in G_N\times G_N:
\atop \{i,j\}\cap\{l,k\}=\emptyset}\left|\hat{\Xi}_N^r\hat{\Lambda}_0^N(i,j,l,k)-\Xi_N^r\Lambda_0^N(i,j,l,k)\right|,
\]
then $h_0^N=0$ according to Assumption (A). By \eqref{equ 4.3}-\eqref{equ 4.10}, for $r\geq 0$,
\begin{align*}
&h_{r+1}^N\\
&\leq \frac{C_{14}}{N}C_{13}Nh_r^N+2C_{15}\frac{C_{14}}{N}C_{12}^r\|\phi\|_\infty+\frac{C_{17}}{N}C_{12}^r\|\phi\|_\infty
+C_{16}h_r^N+2C_{18}\frac{C_{14}}{N}C_{12}^r\|\phi\|_\infty\\
&=C_{19}h_r^N+\frac{C_{20}}{N}C_{12}^r\|\phi\|_\infty,
\end{align*}
where $C_{19}=C_{14}C_{13}+C_{16}$ and $C_{20}=2C_{14}C_{15}+C_{17}+2C_{18}C_{14}$. Then, according to an analysis similar with that leading to \eqref{equ 3.15},
\[
\sup_{(i,j,l,k)\in G_N\times G_N:
\atop \{i,j\}\cap\{l,k\}=\emptyset}\left|\Lambda_t^N(i,j,l,k)-\hat{\Lambda}_t^N(i,j,l,k)\right|\leq \frac{C_{20}\|\phi\|_\infty te^{t\left(C_{12}+C_{19}\right)}}{N}
\]
and hence \eqref{equ 4.11} holds.

\qed

\subsection{Proof of Lemma \ref{lemma 2.3 LLNofVartheta}}\label{subsection 4.2}
In this section we give the proof of Lemma \ref{lemma 2.3 LLNofVartheta}, which follows an analysis similar with that in the proof of Lemma \ref{lemma 2.4 LLNofOmega}.

\proof[Proof of Lemma \ref{lemma 2.3 LLNofVartheta}]
We claim that
\begin{equation}\label{equ 4.14}
\lim_{N\rightarrow+\infty}\sup_{0\leq t\leq T}{\rm Var}\left(\theta_t^N(H_t)\right)=0
\end{equation}
and
\begin{equation}\label{equ 4.15}
\lim_{N\rightarrow+\infty}\sup_{0\leq t\leq T}\left|\mathbb{E}\theta_t^N(H_t)-\theta_t(H_t)\right|=0
\end{equation}
for any $H\in C\left([0, T]\times[0, 1]\right)$. We check \eqref{equ 4.14} and \eqref{equ 4.15} later.
Since
\[
\mathbb{E}\left(\left(\theta_t^N(H_t)-\theta_t(H_t)\right)^2\right)
\leq 2{\rm Var}\left(\theta_t^N(H_t)\right)+2\left(\mathbb{E}\theta_t^N(H_t)-\theta_t(H_t)\right)^2,
\]
Lemma \ref{lemma 2.3 LLNofVartheta} follows from \eqref{equ 4.14} and \eqref{equ 4.15}.

\qed

The remainder of this subsection is devoted to proofs of \eqref{equ 4.14} and \eqref{equ 4.15}.

\proof[Proof of \eqref{equ 4.14}]
According to the definition of $\theta_t^N$,
\[
{\rm Var}\left(\theta_t^N(H_t)\right)\leq \frac{1}{N^2}\sum_{i=1}^N\sum_{j=1}^N\left(\hat{\Lambda}_t^N(i,i,j,j)-\Lambda_t^N(i,i,j,j)\right)H_t\left(\frac{i}{N}\right)H_t\left(\frac{j}{N}\right).
\]
Then by \eqref{equ 4.10},
\[
{\rm Var}\left(\theta_t^N(H_t)\right)\leq \frac{2e^{tC_{12}}}{N}\|H\|_\infty^2+\sup_{(i,j,l,k)\in G_N\times G_N:
\atop \{i,j\}\cap\{l,k\}=\emptyset}\left|\Lambda_t^N(i,j,l,k)-\hat{\Lambda}_t^N(i,j,l,k)\right|\|H\|^2_{\infty}
\]
and hence \eqref{equ 4.14} follows from \eqref{equ 4.11}.

\qed

\proof[Proof of \eqref{equ 4.15}]

According to an analysis similar with that in the proof of Lemma \ref{lemma 2.1 ODE of Hydrodynamic}, the unique solution $\{\theta_t\}_{0\leq t\leq T}$ to Equation \eqref{equ ODE of theta} is given by
\[
\theta_t=e^{tP_2^*}\theta_0+\int_0^te^{(t-s)P_2^*}l_{1,s} ds,
\]
where $\theta_0(du)=\phi(u)du$. Furthermore, $\theta_t(du)=\vartheta(t,u)du$, where $\{\vartheta(t, \cdot)\}_{0\leq t\leq T}$ is the unique solution to the $C[0, 1]$-valued ODE
\[
\begin{cases}
&\frac{d}{dt}\vartheta(t, \cdot)=P_9\vartheta(t, \cdot)+\hat{l}_{1,t},\\
&\vartheta_0=\phi,
\end{cases}
\]
where $\hat{l}_{1,t}\in C[0, 1]$ such that
\begin{align*}
\hat{l}_{1,t}(u)=&2\int_0^1\rho(t,u)\rho(t,v)a_1(u,v)a_2(u,v)\lambda(u,v)dv\\
&+2\int_0^1\lambda(v,u)a_3(v,u)a_4(v,u)\rho(t,u)\rho(t,v)dv
\end{align*}
for all $u\in C[0, 1]$ and $P_9$ is a linear operator from $C[0, 1]$ to $C[0, 1]$ such that
\begin{align*}
P_9f(u)=&\left(c^2(u)-1\right)b(u)f(u)+\int_0^1\lambda(u,v)\left(a_1^2(u,v)-1\right)dvf(u)\\
&\int_0^1\lambda(u,v)a_2^2(u,v)f(v)dv+\int_0^1\lambda(v,u)a_3^2(v,u)f(v)dv\\
&+\int_0^1\lambda(v, u)\left(a_4^2(v,u)-1\right)dvf(u)
\end{align*}
for any $f\in C[0, 1]$ and $u\in [0, 1]$. By \eqref{equ linear system Chapman-Kolmogorov}, for any $i\geq 1$,
\begin{align}\label{equ 4.16}
&\frac{d}{dt}\mathbb{E}\left(\left(X_t^N(i)\right)^2\right)=\\
& b\left(\frac{i}{N}\right)\left(c^2\left(\frac{i}{N}\right)-1\right)\mathbb{E}\left(\left(X_t^N(i)\right)^2\right) \notag\\
&+\frac{1}{N}\sum_{j\neq i}\lambda\left(\frac{i}{N}, \frac{j}{N}\right)\left(a_1^2\left(\frac{i}{N}, \frac{j}{N}\right)-1\right)\mathbb{E}\left(\left(X_t^N(i)\right)^2\right)\notag\\
&+\frac{2}{N}\sum_{j\neq i}\lambda\left(\frac{i}{N}, \frac{j}{N}\right)a_1\left(\frac{i}{N}, \frac{j}{N}\right)a_2\left(\frac{i}{N}, \frac{j}{N}\right)
\mathbb{E}\left(X_t^N(i)X_t^N(j)\right) \notag \\
&+\frac{1}{N}\sum_{j\neq i}\lambda\left(\frac{i}{N}, \frac{j}{N}\right)a_2^2\left(\frac{i}{N}, \frac{j}{N}\right)\mathbb{E}\left(\left(X_t^2(j)\right)^2\right) \notag\\
&+\frac{1}{N}\sum_{j\neq i}\lambda\left(\frac{j}{N}, \frac{i}{N}\right)\left(a_4^2\left(\frac{j}{N}, \frac{i}{N}\right)-1\right)\mathbb{E}\left(\left(X_t^N(i)\right)^2\right) \notag\\
&+\frac{2}{N}\sum_{j\neq i}\lambda\left(\frac{j}{N}, \frac{i}{N}\right)a_3\left(\frac{j}{N}, \frac{i}{N}\right)a_4\left(\frac{j}{N}, \frac{i}{N}\right)
\mathbb{E}\left(X_t^N(i)X_t^N(j)\right) \notag\\
&+\frac{1}{N}\sum_{j\neq i}\lambda\left(\frac{j}{N}, \frac{i}{N}\right)a_3^2\left(\frac{j}{N}, \frac{i}{N}\right)\mathbb{E}\left(\left(X_t^2(j)\right)^2\right). \notag
\end{align}
Let $\hat{\theta}_t^N$ be the element in $\left(C[0, 1]\right)^\prime$ such that $\hat{\theta}_t^N(f)=\mathbb{E}\theta_t^N(f)$ for any $f\in C[0, 1]$, then
by \eqref{equ 3.8} and \eqref{equ 4.16}, $\left\{\hat{\theta}_t^N(f):~0\leq t\leq T\right\}_{N\geq 1}$ are uniformly bounded and equicontinuous for any $f\in C[0, 1]$. Then, by Ascoli-Arzela theorem, any subsequence $\{\hat{\theta}_t^{N_k}:~0\leq t\leq T\}_{k\geq 1}$ of $\{\hat{\theta}_t^N:~0\leq t\leq T\}_{N\geq 1}$ has a subsequence $\{\hat{\theta}_t^{N_{k_j}}:~0\leq t\leq T\}_{j\geq 1}$ such that $\hat{\theta}_t^{N_{k_j}}(f)$ converges to $\tau_t(f)$ uniformly in $t\in[0, T]$ as $j\rightarrow+\infty$ for some $\{\tau_t\}_{0\leq t\leq T}\in \mathcal{D}\left([0, T], \left(C[0, 1]\right)^\prime\right)$ and any $f\in C[0, 1]$. As a result, we only need to show that $\tau=\theta$ to complete this proof. By \eqref{equ 3.8}, \eqref{equ 4.16} and Assumption (A),
\begin{align}\label{equ 4.17}
\hat{\theta}_t^{N_{k_j}}(f)=&\frac{1}{N_{k_j}}\sum_{i=1}^{N_{k_j}}\phi\left(\frac{i}{N_{k_j}}\right)f\left(\frac{i}{N_{k_j}}\right)
+\int_0^t\hat{\theta}_s^{N_{k_j}}(P_{2,N_{k_j}}f)ds \notag\\
&+\int_0^t\hat{\omega}_s^{N_{k_j}}(\mathcal{K}f)ds+O\left(\frac{1}{N_{k_j}}\right)
\end{align}
for any $0\leq t\leq T$ and $f\in C[0, 1]$, where $\mathcal{K}$ is the linear operator from $C[0, 1]$ to $C\left([0, 1]\times [0, 1]\right)$ such that
\[
\mathcal{K}f(u, v)=2\lambda(u,v)a_1(u,v)a_2(u,v)f(u)+2\lambda(v,u)a_3(v,u)a_4(v,u)f(u)
\]
for any $f\in C[0, 1], u,v\in [0, 1]$ and $P_{2,N}$ is the linear operator from $C[0, 1]$ to $C[0, 1]$ such that
\begin{align*}
&P_{2,N}f(u)=\\
&b(u)\left(c^2(u)-1\right)f(u)+\frac{f(u)}{N}\sum_{j=1}^N\lambda\left(u, \frac{j}{N}\right)\left(a_1^2\left(u, \frac{j}{N}\right)-1\right)\\
&+\frac{1}{N}\sum_{j=1}^N\lambda\left(\frac{j}{N}, u\right)a_2^2\left(\frac{j}{N}, u\right)f\left(\frac{j}{N}\right)\\
&+\frac{1}{N}\sum_{j=1}^N\lambda\left(u, \frac{j}{N}\right)a_3^2\left(u, \frac{j}{N}\right)f\left(\frac{j}{N}\right)\\
&+\frac{f(u)}{N}\sum_{j=1}^N\lambda\left(\frac{j}{N}, u\right)\left(a_4^2\left(\frac{j}{N}, u\right)-1\right)
\end{align*}
for any $f\in C[0, 1], u\in [0, 1]$. According to definitions of $P_2, P_{2,N}$,
\[
 \lim_{N\rightarrow+\infty}\left\|P_2f-P_{2,N}f\right\|_\infty=0
\]
for any $f\in C[0, 1]$. By Lemma \ref{lemma 2.4 LLNofOmega},
\[
\lim_{N\rightarrow+\infty}\hat{\omega}_s^N(\mathcal{K}f)=\int_0^1\int_0^1\rho(s,u)\rho(s,v)\mathcal{K}f(u,v)dudv=l_{1,s}(f)
\]
uniformly for $s\in [0, T]$. Then let $j\rightarrow+\infty$ in \eqref{equ 4.17}, by \eqref{equ 3.8},
\begin{equation}\label{equ 4.18}
\tau_t(f)=\int_0^1\phi(u)f(u)du+\int_0^t\tau_s(P_2f)ds+\int_0^tl_{1,s}(f)ds
\end{equation}
for any $0\leq t\leq T$ and $f\in C[0, 1]$. Since $\theta$ is the solution to Equation \eqref{equ ODE of theta}, \eqref{equ 4.18} still holds when we replace $\tau$ by $\theta$ and hence $\tau=\theta$ according to Gronwall's inequality.

\qed

\subsection{Proof of Theorem \ref{theorem 2.6 fluctuation of the NurnLinearSyetem}}\label{subsection 4.3}
In this subsection we give the proof of Theorem \ref{theorem 2.6 fluctuation of the NurnLinearSyetem}. As a preliminary, we first prove Lemma \ref{lemma 2.5}.

\proof[Proof of Lemma \ref{lemma 2.5}]
We only need to show that $\left[f,f\right]_s\geq 0$ for any $f\in C[0, 1]$. For any $f\in C[0, 1], 0\leq s\leq T$ and $N\geq 1$, let
\begin{align*}
&\mathcal{Z}_s^N(f)=\\
&\sum_{i=1}^Nb\left(\frac{i}{N}\right)\left(\frac{f\left(\frac{i}{N}\right)\left(c\left(\frac{i}{N}\right)-1\right)}{\sqrt{N}}X_s^N(i)\right)^2\\
&+\frac{1}{N}\sum_{i=1}^N\sum_{j\neq i}\lambda\left(\frac{i}{N}, \frac{j}{N}\right)\Bigg(\\
&\text{\quad}\frac{f\left(\frac{i}{N}\right)}{\sqrt{N}}\left(a_1\left(\frac{i}{N}, \frac{j}{N}\right)X_s^N(i)+a_2\left(\frac{i}{N}, \frac{j}{N}\right)X_s^N(j)-X_s^N(i)\right) \\
&\text{\quad\quad}+\frac{f\left(\frac{j}{N}\right)}{\sqrt{N}}\left(a_3\left(\frac{i}{N}, \frac{j}{N}\right)X_s^N(i)+a_4\left(\frac{i}{N}, \frac{j}{N}\right)X_s^N(j)-X_s^N(j)\right)\Bigg)^2,
\end{align*}
then $\mathcal{Z}_s^N(f)\geq 0$. By direct calculation and \eqref{equ 4.10},
\begin{equation}\label{equ 4.19}
\mathcal{Z}_s^N(f)=\theta_t^N\left(P_{10}f+P_{11, N}f\right)+\omega_t^N(P_{12}f)+\zeta_s^N(f),
\end{equation}
where
\[
\sup_{0\leq s\leq T}\mathbb{E}\left(\left(\zeta_s^N(f)\right)^2\right)=O\left(\frac{1}{N}\right)
\]
and $P_{10}, P_{11, N}$ are nonlinear operators from $C[0, 1]$ to $C[0, 1]$ such that
\[
P_{10}f(u)=b(u)\left(c(u)-1\right)^2f^2(u)
\]
and
\begin{align*}
&P_{11, N}f(u)=\\
&\frac{f^2(u)}{N}\sum_{j=1}^N\lambda\left(u, \frac{j}{N}\right)\left(a_1\left(u, \frac{j}{N}\right)-1\right)^2+\frac{1}{N}\sum_{j=1}^N\lambda\left(\frac{j}{N}, u\right)a_2^2\left(\frac{j}{N}, u\right)f^2\left(\frac{j}{N}\right)\\
&+\frac{1}{N}\sum_{j=1}^N\lambda\left(u, \frac{j}{N}\right)a_3^2\left(u, \frac{j}{N}\right)f^2\left(\frac{j}{N}\right)+\frac{f^2(u)}{N}\sum_{j=1}^N\lambda\left(\frac{j}{N},u\right)\left(a_4\left(\frac{j}{N}, u\right)-1\right)^2\\
&+\frac{2}{N}\sum_{j=1}^N\lambda\left(u, \frac{j}{N}\right)\left(a_1\left(u,\frac{j}{N}\right)-1\right)a_3\left(u, \frac{j}{N}\right)f(u)f\left(\frac{j}{N}\right)\\
&+\frac{2}{N}\sum_{j=1}^N\lambda\left(\frac{j}{N}, u\right)a_2\left(\frac{j}{N},u\right)\left(a_4\left(\frac{j}{N},u\right)-1\right)f(u)f\left(\frac{j}{N}\right)
\end{align*}
for any $u\in [0, 1], f\in C[0, 1]$ and $P_{12}$ is the nonlinear operator from $C[0, 1]$ to $C\left([0, 1]\times[0, 1]\right)$ such that
\begin{align*}
P_{12}f(u,v)=&2\lambda(u,v)\left(a_1(u,v)-1\right)a_2(u,v)f^2(u)\\
&+2\lambda(u,v)a_3(u,v)\left(a_4(u,v)-1\right)f^2(v)\\
&+2\lambda(u,v)\left(a_1(u,v)-1\right)\left(a_4(u,v)-1\right)f(u)f(v)\\
&+2\lambda(u,v)a_2(u,v)a_3(u,v)f(u)f(v)
\end{align*}
for any $f\in C[0, 1], u, v\in [0, 1]$. According to the definition of $P_{11, N}$, for any $f\in C[0, 1]$,
\[
\lim_{N\rightarrow+\infty}\left\|P_{11, N}f-P_{11}f\right\|_\infty=0,
\]
where $P_{11}$ is the nonlinear operator from $C[0, 1]$ to $C[0, 1]$ such that
\begin{align*}
&P_{11}f(u)=\\
&f^2(u)\int_0^1\lambda(u,v)\left(a_1(u,v)-1\right)^2dv+\int_0^1\lambda(v,u)a_2^2(v,u)f^2(v)dv\\
&+\int_0^1\lambda(u,v)a_3^2(u,v)f^2(v)dv+f^2(u)\int_0^1\lambda(v,u)\left(a_4(v,u)-1\right)^2dv\\
&+2f(u)\int_0^1\lambda(u,v)\left(a_1(u,v)-1\right)a_3(u,v)f(v)dv\\
&+2f(u)\int_0^1\lambda(v,u)a_2(v,u)\left(a_4(v,u)-1\right)f(v)dv
\end{align*}
for any $f\in C[0, 1], u\in [0, 1]$. Then, by \eqref{equ 4.10}, \eqref{equ 4.19} and Lemmas \ref{lemma 2.3 LLNofVartheta}, \ref{lemma 2.4 LLNofOmega}
\begin{align}\label{equ 4.20}
\lim_{N\rightarrow+\infty}\sup_{0\leq s\leq T}\mathbb{E}\Bigg(\Big(&\mathcal{Z}_s^N(f)-\theta_s(P_{10}f+P_{11}f)\\
&-\int_0^1\int_0^1\rho(s,u)\rho(s,v)P_{12}f(u,v)dudv\Big)^2\Bigg)=0. \notag
\end{align}
By direct calculation,
\[
\left[f,f\right]_s=\theta_s(P_{10}f+P_{11}f)+\int_0^1\int_0^1\rho(s,u)\rho(s,v)P_{12}f(u,v)dudv
\]
and then $\left[f,f\right]_s\geq 0$ follows from \eqref{equ 4.20} and the fact that $\mathcal{Z}_s^N(f)\geq 0$.

\qed

To prove Theorem \ref{theorem 2.6 fluctuation of the NurnLinearSyetem}, we need the following lemma.

\begin{lemma}\label{lemma 4.1}
Under Assumption (A), $\{V_t^N:~0\leq t\leq T\}_{N\geq 1}$ are tight.
\end{lemma}

We prove Lemma \ref{lemma 4.1} at the end of this subsection. Now we give the proof of Theorem \ref{theorem 2.6 fluctuation of the NurnLinearSyetem}.

\proof[Proof of Theorem \ref{theorem 2.6 fluctuation of the NurnLinearSyetem}]

By Lemma \ref{lemma 4.1}, let $\hat{V}=\{\hat{V}_t\}_{0\leq t\leq T}$ be a weak limit of a subsequence $\{V_t^{N_j}:~0\leq t\leq T\}_{j\geq 1}$ of $\{V_t^N:~0\leq t\leq T\}_{N\geq 1}$. We only need to show that $\hat{V}$ satisfies all the three properties in the definition of the solution to the martingale problem related to the generalized O-U process in \eqref{equ 2.6 generalized O-U}. The third property of the initial state $\hat{V}_0$ follows directly from Assumption (A), so we only need to check the first and the second properties. For any $f\in C[0, 1]$, the first property of $\hat{V}$ that $\{\hat{V}_t(f)\}_{0\leq t\leq T}$ is continuous follows from
\begin{equation}\label{equ 4.21}
\lim_{N\rightarrow+\infty}\sup_{0\leq t\leq T}\left|V_t^N(f)-V_{t-}^N(f)\right|=0
\end{equation}
in probability. Now we check \eqref{equ 4.21}. According to transition rates of $\{X_t^N\}_{t\geq 0}$,
\[
\sup_{0\leq t\leq T}\left|V_t^N(f)-V_{t-}^N(f)\right|\leq C_{21}\frac{\|f\|_\infty}{\sqrt{N}}\sup_{0\leq t\leq T,\atop 1\leq i\leq N}X_t^N(i),
\]
where $C_{21}=\|C\|_\infty+\sum_{k=1}^4\|a_k\|_\infty+3$. So we only need to check
\begin{equation}\label{equ 4.22}
\lim_{N\rightarrow+\infty}\frac{1}{\sqrt{N}}\sup_{1\leq i\leq N,\atop 0\leq t\leq T}X_t^N(i)=0
\end{equation}
in probability to prove \eqref{equ 4.21}. Let $\{\hat{X}_t^N\}_{t\geq 0}$ be our process with parameters
\[
\hat{b}, \hat{c}, \hat{\lambda}, \hat{a}_1, \hat{a}_2, \hat{a}_3, \hat{a}_4
\]
given by
\[
\hat{b}=\|b\|_\infty, \hat{\lambda}=\|\lambda\|_\infty, \hat{c}=\|c\|_\infty+1
\]
and
\[
\hat{a}_1=\|a_1\|_\infty+1, \hat{a}_2=\|a_2\|_\infty+1, \hat{a}_3=\|a_3\|_\infty+1, \hat{a}_4=\|a_4\|_\infty+1.
\]
Then $X_t^N(i)\leq \hat{X}_t^N(i)$ for all $t\geq 0$ in the sense of coupling and $\hat{X}_t^N(i)$ is increasing with $t$. Hence,
\begin{equation}\label{equ 4.23}
\left(\frac{1}{\sqrt{N}}\sup_{1\leq i\leq N,\atop 0\leq t\leq T}X_t^N(i)\right)^4\leq \frac{1}{N^2}\sum_{i=1}^N\left(\hat{X}_T^N(i)\right)^4.
\end{equation}
By \eqref{equ 4.10},
\[
\sup_{1\leq i\leq N}\mathbb{E}\left(\left(\hat{X}_T^N(i)\right)^4\right)<e^{T\hat{C}_{12}},
\]
where $\hat{C}_{12}=C_{12}\left(\hat{b},\hat{c},\hat{\lambda}, \hat{a_1}, \hat{a_2}, \hat{a_3}, \hat{a_4}\right)$. Then by \eqref{equ 4.23}, $\left(\frac{1}{\sqrt{N}}\sup_{1\leq i\leq N,\atop 0\leq t\leq T}X_t^N(i)\right)^4$ converges to $0$ in $L^1$ and hence \eqref{equ 4.22} holds. In conclusion, $\{\hat{V}_t(f)\}_{0\leq t\leq T}$ is continuous.

\quad

At last we check the second property of $\hat{V}$, i.e.,
\begin{align*}
\Bigg\{G\left(\hat{V}_t(f)\right)-G\left(\hat{V}_0(f)\right)&-\int_0^tG^\prime\left(\hat{V}_s(f)\right)\hat{V}_s\left(P_1f\right)ds \\
&-\frac{1}{2}\int_0^tG^{\prime\prime}\left(\hat{V}_s(f)\right)[f,f]_sds\Bigg\}_{0\leq t\leq T}
\end{align*}
is a martingale for any $f\in C[0, 1]$ and $G\in C_c^{\infty}(\mathbb{R})$. For $f\in C[0, 1]$. let
\begin{align*}
\Upsilon_t^N(f)=G\left(V_t^N(f)\right)-G\left(V_0^N(f)\right)&-\int_0^t\left(\partial_s+\mathcal{L}_N\right)G(V_s^N(f))ds,
\end{align*}
then $\{\Upsilon_t^N(f)\}_{0\leq t\leq T}$ is a martingale according to Dynkin's martingale formula. By the definition of $\mathcal{L}_N$, \eqref{equ 3.8} and Taylor's expansion up to the second order,
\begin{equation}\label{equ 4.24}
\left(\partial_s+\mathcal{L}_N\right)G\left(V_s^N(f)\right)
=G^\prime\left(V_s^N(f)\right)V_s^N\left(P_{1,N}f\right)+\frac{G^{\prime\prime}\left(V_s^N(f)\right)}{2}\mathcal{Z}_s^N(f)+\varsigma_s^N+\varphi_s^N,
\end{equation}
where
\begin{align*}
&|\varsigma_s^N|\leq \frac{\left\|G^{\prime\prime\prime}\right\|_\infty}{6}\Bigg(\frac{\|b\|_\infty\|f\|_\infty\left(\|c\|_\infty+1\right)}
{N^{\frac{3}{2}}}\sum_{i=1}^N\left(X_s^N(i)\right)^3\\
&+\frac{\|\lambda\|_\infty\|f\|_\infty}{N^{\frac{5}{2}}}\sum_{i,j}\left((\|a_1\|_\infty+\|a_3\|_\infty+1)X_t^N(i)
+(\|a_2\|_\infty+\|a_4\|_\infty+1)X_t^N(j)\right)^3\Bigg)
\end{align*}
and
\[
\sup_{0\leq s\leq T}\mathbb{E}|\varphi_s^N|=O\left(\frac{1}{\sqrt{N}}\right).
\]
Then by \eqref{equ 4.10}, for $0\leq t\leq T$,
\[
\lim_{N\rightarrow+\infty}\int_0^t(\varsigma_s^N+\varphi_s^N)ds=0
\]
in $L^1$. By \eqref{equ 3.15} and Cauchy-Schwarz inequality, for $0\leq t\leq T$,
\[
\lim_{N\rightarrow+\infty}\left|\int_0^t G^\prime\left(V_s^N(f)\right)\left(V_s^N\left(P_{1,N}f\right)-V_s^N(P_1f)\right)ds\right|=0
\]
in $L^2$.
By \eqref{equ 4.20} and Cauchy-Schwarz inequality, for $0\leq t\leq T$,
\[
\lim_{N\rightarrow+\infty}\left|\int_0^t\frac{1}{2}G^{\prime\prime}\left(V_s^N(f)\right)\left(\mathcal{Z}_s^N(f)-\left[f,f\right]_s\right)ds\right|=0
\]
in $L^2$. As a result, let $N\rightarrow+\infty$ in \eqref{equ 4.24}, then $\Upsilon_t^{N_j}(f)$ converges weakly to $\hat{\Upsilon}_t(f)$ as $j\rightarrow+\infty$ and $\{\Upsilon_t^{N_j}(f)\}_{j\geq 1}$ are uniformly integrable for any $f\in C[0, 1], 0\leq t\leq T$, where
\begin{align*}
\hat{\Upsilon}_t(f)=G\left(\hat{V}_t(f)\right)-G\left(\hat{V}_0(f)\right)&-\int_0^tG^\prime\left(\hat{V}_s(f)\right)\hat{V}_s\left(P_1f\right)ds \\
&-\frac{1}{2}\int_0^tG^{\prime\prime}\left(\hat{V}_s(f)\right)[f,f]_sds.
\end{align*}
As a result, by Theorem 5.3 of \cite{Whitt2007}, $\{\hat{\Upsilon}_t(f)\}_{0\leq t\leq T}$ is a martingale for any $f\in C[0, 1]$ and the proof is complete.

\qed

At last, we prove Lemma \ref{lemma 4.1}.

\proof[Proof of Lemma \ref{lemma 4.1}]

By Aldous' criteria, we only need to check that
\begin{equation}\label{equ 4.25}
\lim_{M\rightarrow+\infty}\limsup_{N\rightarrow+\infty}P\left(\left|V_t^N(f)\right|\geq M\right)=0
\end{equation}
for any $0\leq t\leq T, f\in C[0, 1]$ and
\begin{equation}\label{equ 4.26}
\lim_{\delta\rightarrow0}\limsup_{N\rightarrow+\infty}\sup_{\upsilon\in\mathcal{T},s<\delta}
P\left(\left|V_{\upsilon+s}^N(f)-V_\upsilon^N(f)\right|>\epsilon\right)=0
\end{equation}
for any $\epsilon>0, f\in C[0, 1]$, where $\mathcal{T}$ is the set of stopping times of $\{X_t^N\}_{t\geq 0}$ bounded by $T$. By Dynkin's martingale formula,
\[
V_t^N(f)=V_0^N(f)+\eta_t^N(f)+\sigma_t^N(f),
\]
where
\[
\eta_t^N(f)=\int_0^t\left(\partial_s+\mathcal{L}_N\right)V_s^N(f)ds
\]
and $\{\sigma_t^N(f)\}_{0\leq t\leq T}$ is a martingale with quadratic variation process $\langle \sigma^N(f)\rangle_t$ given by
\[
\langle\sigma^N(f)\rangle_t=\int_0^t\left(\mathcal{L}_N\left(\left(V_s^N(f)\right)^2\right)-2V_s^N(f)\mathcal{L}_NV_s^N(f)\right)ds.
\]
Hence, to check \eqref{equ 4.25} and \eqref{equ 4.26}, we only need to show that the following four equations hold.

1) For any $0\leq t\leq T$ and $f\in C[0, 1]$,
\begin{equation}\label{equ 4.27}
\lim_{M\rightarrow+\infty}\limsup_{N\rightarrow+\infty}P\left(\left|\eta_t^N(f)\right|\geq M\right)=0.
\end{equation}

2) For any $\epsilon>0$ and $f\in C[0, 1]$,
\begin{equation}\label{equ 4.28}
\lim_{\delta\rightarrow0}\limsup_{N\rightarrow+\infty}\sup_{\upsilon\in\mathcal{T},s<\delta}
P\left(\left|\eta_{\upsilon+s}^N(f)-\eta_\upsilon^N(f)\right|>\epsilon\right)=0.
\end{equation}

3) For any $0\leq t\leq T$ and $f\in C[0, 1]$,
\begin{equation}\label{equ 4.29}
\lim_{M\rightarrow+\infty}\limsup_{N\rightarrow+\infty}P\left(\left|\sigma_t^N(f)\right|\geq M\right)=0.
\end{equation}

4) For any $\epsilon>0$ and $f\in C[0, 1]$,
\begin{equation}\label{equ 4.30}
\lim_{\delta\rightarrow0}\limsup_{N\rightarrow+\infty}\sup_{\upsilon\in\mathcal{T},s<\delta}
P\left(\left|\sigma_{\upsilon+s}^N(f)-\sigma_\upsilon^N(f)\right|>\epsilon\right)=0.
\end{equation}

Now we check \eqref{equ 4.27}-\eqref{equ 4.30}. According to the definition of $\mathcal{L}_N$ and \eqref{equ linear system Chapman-Kolmogorov},
\[
\eta_t^N(f)=\int_0^t\frac{1}{\sqrt{N}}\sum_{i=1}^N\left(X_s^N(i)-\mathbb{E}X_s^N(i)\right)\left(H_N^f(i)\right)ds,
\]
where
\begin{align*}
&H_N^f(i)=\\
&\frac{f(i)}{N}\sum_{j\neq i}\lambda\left(\frac{i}{N}, \frac{j}{N}\right)\left(a_1\left(\frac{i}{N}, \frac{j}{N}\right)-1\right)+b\left(\frac{i}{N}\right)\left(c\left(\frac{i}{N}\right)-1\right)f\left(\frac{i}{N}\right)\\
&+\frac{1}{N}\sum_{j\neq i}\lambda\left(\frac{i}{N}, \frac{j}{N}\right)a_3\left(\frac{i}{N}, \frac{j}{N}\right)f\left(\frac{j}{N}\right)\\
&+\frac{f\left(\frac{i}{N}\right)}{N}\sum_{j\neq i}\lambda\left(\frac{j}{N}, \frac{i}{N}\right)\left(a_4\left(\frac{j}{N}, \frac{i}{N}\right)-1\right)\\
&+\frac{1}{N}\sum_{j\neq i}\lambda\left(\frac{j}{N}, \frac{i}{N}\right)f\left(\frac{j}{N}\right)a_2\left(\frac{j}{N}, \frac{i}{N}\right)
\end{align*}
for all $1\leq i\leq N$. Since
\[
\sup_{1\leq i\leq N}\left|H_N^f(i)\right|\leq C_{22}\|f\|_\infty
\]
for some $C_{22}<+\infty$ independent of $N, f$, by \eqref{equ 3.8}, \eqref{equ 3.15} and Cauchy-Schwarz inequality,
\begin{align*}
&\mathbb{E}\left(\left(\eta_t^N(f)\right)^2\right)\\
&\leq \frac{T}{N}\int_0^t\sum_{i=1}^N\sum_{j=1}^N{\rm Cov}\left(X_s^N(i), X_s^N(j)\right)H_N^f(i)H_N^f(j)ds\\
&\leq \frac{T^2C_{22}^2\|f\|_\infty^2}{N}\left(Ne^{C_5T}\|\phi\|_\infty+N^2\frac{C_{11}\|\phi\|_\infty Te^{(C_7+C_5)T}}{N}\right)\\
&=\left(T^2C_{22}^2\|f\|_\infty^2\right)\left(e^{C_5T}\|\phi\|_\infty+C_{11}\|\phi\|_\infty Te^{(C_7+C_5)T}\right).
\end{align*}
Then, \eqref{equ 4.27} follows from Markov's inequality.

\quad

Similarly, by Cauchy-Schwarz inequality,
\begin{align*}
&\mathbb{E}\left(\left|\eta_{\upsilon+s}^N(f)-\eta_\upsilon^N(f)\right|^2\right)\\
&=\mathbb{E}\left(\left(\int_0^{T+\delta}1_{\{\upsilon\leq u\leq \upsilon+s\}}\frac{1}{\sqrt{N}}\sum_{i=1}^N\left(X_u^N(i)-\mathbb{E}X_u^N(i)\right)H^f_N(i)du\right)^2\right)\\
&\leq \mathbb{E}\left(\int_0^{T+\delta}1_{\{\upsilon\leq u\leq \upsilon+s\}}du\int_0^{T+\delta}\left(\frac{1}{\sqrt{N}}\sum_{i=1}^N\left(X_u^N(i)-\mathbb{E}X_u^N(i)\right)H^f_N(i)\right)^2du\right)\\
&=\mathbb{E}\left(s\int_0^{T+\delta}\left(\frac{1}{\sqrt{N}}\sum_{i=1}^N\left(X_u^N(i)-\mathbb{E}X_u^N(i)\right)H^f_N(i)\right)^2du\right)\\
&=s\left(\int_0^{T+\delta}\frac{1}{N}\sum_{i=1}^N\sum_{j=1}^N{\rm Cov}\left(X_u^N(i), X_u^N(j)\right)H_N^f(i)H_N^f(j)du\right)\\
&\leq \delta \left((T+\delta)C_{22}^2\|f\|_\infty^2\right)\left(e^{C_5(T+\delta)}\|\phi\|_\infty+C_{11}\|\phi\|_\infty Te^{(C_7+C_5)(T+\delta)}\right)
\end{align*}
for $0\leq s\leq \delta$. As result, \eqref{equ 4.28} follows from Markov's inequality.

\quad

By direct calculation,
\begin{align*}
\langle \sigma^N(f)\rangle_t&=\int_0^t\frac{1}{N^2}\sum_{i=1}^N\sum_{j\neq i}\lambda\left(\frac{i}{N}, \frac{j}{N}\right)\times \\
&\Bigg(f\left(\frac{i}{N}\right)\left(\left(a_1\left(\frac{i}{N}, \frac{j}{N}\right)-1\right)X_s^N(i)+a_2\left(\frac{i}{N}, \frac{j}{N}\right)X_s^N(j)\right)\\
&+f\left(\frac{j}{N}\right)\left(\left(a_4\left(\frac{i}{N}, \frac{j}{N}\right)-1\right)X_s^N(j)+a_3\left(\frac{i}{N}, \frac{j}{N}\right)X_s^N(i)\right)\Bigg)^2 ds.
\end{align*}
As a result, there exists $C_{23}<+\infty$ independent of $N, t$ such that
\[
\langle \sigma^N(f)\rangle_t-\langle\sigma^N(f)\rangle_s\leq C_{23}\|f\|_\infty^2\int_s^t \theta_u^N(\vec{1})du
\]
for any $s\leq t$, where $\vec{1}$ is the constant function taking value $1$. Then, by \eqref{equ 3.8},
\[
\mathbb{E}\left(\left(\sigma_t^N(f)\right)^2\right)=\mathbb{E}\langle\sigma^N(f)\rangle_t\leq C_{23}\|f\|_\infty^2Te^{C_5T}\|\phi\|_\infty
\]
and hence \eqref{equ 4.29} follows from Markov's inequality.

\quad

By Cauchy-Schwarz inequality and \eqref{equ 4.10}, for $0\leq s\leq \delta$,
\begin{align*}
&\mathbb{E}\left(\left(\sigma_{\upsilon+s}^N(f)-\sigma_\upsilon^N(f)\right)^2\right)=\mathbb{E}\left(\langle \sigma^N(f)\rangle_{\upsilon+s}-\langle \sigma^N(f)\rangle_\upsilon\right)\\
&\leq C_{23}\|f\|_\infty^2\mathbb{E}\left(\int_\upsilon^{\upsilon+s}\theta_u^N(\vec{1})du\right)\leq C_{23}\|f\|_\infty^2\mathbb{E}\left(\sqrt{s}\sqrt{\int_0^{T+\delta}\left(\theta_u^N(\vec{1})\right)^2du}\right)\\
&\leq C_{23}\|f\|_\infty^2\sqrt{\delta}\sqrt{\int_0^{T+\delta}\mathbb{E}\left(\left(\theta_u^N(\vec{1})\right)^2\right)du}\\
&\leq C_{23}\|f\|_\infty^2\sqrt{\delta}\sqrt{(T+\delta)e^{(T+\delta)C_{12}}}.
\end{align*}
As a result, \eqref{equ 4.30} follows from Markov's inequality. Since \eqref{equ 4.27}-\eqref{equ 4.30} all hold, \eqref{equ 4.25} and \eqref{equ 4.26} hold and the proof is complete.

\qed

\section{Applications}\label{section five}
In this section we apply our main results in three Examples given in Section \ref{section one}, i.e., $N$-urn voter models, pair-symmetric $N$-urn exclusion processes and $N$-urn binary contact path processes.

\quad

\textbf{Example 1} \emph{$N$-urn voter model}. In this case $b=a_1=a_3=0$ and $a_2=a_4=1$. Then, $P_1$ is given by
\[
P_1f(u)=\int_0^1 \lambda(v,u)f(v)dv-\int_0^1\lambda(u,v)dvf(u)
\]
for any $u\in [0, 1]$ and $f\in C[0, 1]$. Furthermore, by direct calculation, $P_2=P_1$, $l_1=0$, $\{\vartheta(t,\cdot)\}_{0\leq t\leq T}=\{\rho(t,\cdot)\}_{0\leq t\leq T}$ and
\[
\left[f, f\right]_s=\int_0^1\int_0^1f^2(u)\left[\rho(s,u)\lambda(u,v)+\rho(s,v)\lambda(u,v)-2\lambda(u,v)\rho(s,v)\rho(s,u)\right]dudv
\]
for any $f\in C[0, 1]$ and $0\leq s\leq T$. As a result, by Theorem \ref{theorem 2.2 hydrodynamic limit},
\[
\lim_{N\rightarrow+\infty}\mu_t^N(f)=\int_0^1\rho(t,u)f(u)du
\]
in $L^2$ for any $t\geq 0$ and $f\in C[0, 1]$,  where $\{\rho(t,\cdot)\}_{t\geq 0}$ is given by
\[
\begin{cases}
&\frac{d}{dt}\rho(t, \cdot)=\int_0^1\lambda(\cdot, v)\rho(t, v)dv-\int_0^1\lambda(\cdot, v)dv\rho(t, \cdot),\\
&\rho(0, \cdot)=\phi.
\end{cases}
\]
By Theorem \ref{theorem 2.6 fluctuation of the NurnLinearSyetem}, $V^N$ converges weakly to $V$ as $N\rightarrow+\infty$,
where $V=\{V_t\}_{0\leq t\leq T}$ is the unique solution to
\[
\begin{cases}
&dV_t=P_1^{*}V_tdt+\mathcal{R}_t^*d\mathcal{W}_t\text{~for~}0\leq t\leq T,\\
&V_0(f) \text{~follows~}\mathbb{N}\left(0, \int_0^1f^2(u)\phi(u)(1-\phi(u))du\right)\text{~for any~}f\in C[0, 1],\\
&V_0 \text{~is independent of~}\{C_t\}_{0\leq t\leq T},
\end{cases}
\]
where $\{\mathcal{W}_t\}_{t\geq 0}$ is the standard $\left(C[0, 1]\right)^\prime$-valued Brownian motion such that $\{\mathcal{W}_t(f)\}_{t\geq 0}$ is a real-valued Brownian motion with $\mathcal{W}_0(f)=0, {\rm Var}\left(\mathcal{W}_t(f)\right)=t\int_0^1f^2(u)du$ for any $f\in C[0, 1], t\geq 0$ and
$\mathcal{R}_t$ is a linear operator from $C[0, 1]$ to $C[0, 1]$ such that
\[
\mathcal{R}_tf(u)=\sqrt{\int_0^1\rho(t,u)\lambda(u,v)+\rho(t,v)\lambda(u,v)-2\lambda(u,v)\rho(t,v)\rho(t,u)dv}f(u)
\]
for any $f\in C[0, 1], u\in [0, 1], t\geq 0$. Especially, when $\lambda=\lambda_0\in (0, +\infty)$ and $\phi=p\in [0, 1]$, we have $\rho(t,\cdot)\equiv p$ for any $t\geq 0$ and $\mathcal{R}_t^*d\mathcal{W}_t=\sqrt{2\lambda_0p(1-p)}d\mathcal{W}_t$.

\quad

\textbf{Example 2} \emph{Pair-symmetric $N$-urn exclusion process}. In this case, $b=a_1=a_4=0$ and $a_2=a_3=1$. Hence, $P_1$ is given by
\begin{align*}
P_1f(u)=&-f(u)\int_0^1\lambda(u, v)dv+\int_0^1\lambda(u,v)f(v)dv\\
&-f(u)\int_0^1\lambda(v,u)dv+\int_0^1\lambda(v,u)f(v)dv.
\end{align*}
Furthermore, by direct calculation, $P_2=P_1$, $l_1=0$, $\{\vartheta(t,\cdot)\}_{0\leq t\leq T}=\{\rho(t,\cdot)\}_{0\leq t\leq T}$ and
\[
\left[f, f\right]_s=\int_0^1\int_0^1\lambda(u,v)\left(\rho(t,u)+\rho(t,v)-2\rho(t,u)\rho(t,v)\right)\left(f(u)-f(v)\right)^2dudv
\]
for any $f\in C[0, 1]$ and $0\leq s\leq T$. As a result, by Theorem \ref{theorem 2.2 hydrodynamic limit},
\[
\lim_{N\rightarrow+\infty}\mu_t^N(f)=\int_0^1f(u)\rho(t,u)du
\]
in $L^2$ for any $f\in C[0, 1]$ and $t\geq 0$, where $\rho$ is given by
\[
\begin{cases}
&\frac{d}{dt}\rho(t,\cdot)=-\int_0^1\lambda(v,\cdot)dv\rho(t,\cdot)-\rho(t,\cdot)\int_0^1\lambda(\cdot, v)dv\\
&\text{\quad\quad\quad\quad~}+\int_0^1\lambda(\cdot, v)\rho(t,v)dv+\int_0^1\lambda(v,\cdot)\rho(t,v)dv,\\
&\rho_0=\phi.
\end{cases}
\]

By Theorem \ref{theorem 2.6 fluctuation of the NurnLinearSyetem}, $V^N$ converges weakly to $V$ as $N\rightarrow+\infty$, where $V$ is the solution to
\[
\begin{cases}
&dV_t=P_1^{*}V_tdt+\mathcal{U}_t^*d\mathcal{B}_t\text{~for~}0\leq t\leq T,\\
&V_0(f) \text{~follows~}\mathbb{N}\left(0, \int_0^1f^2(u)\phi(u)(1-\phi(u))du\right)\text{~for any~}f\in C[0, 1],\\
&V_0 \text{~is independent of~}\{C_t\}_{0\leq t\leq T},
\end{cases}
\]
where $\{\mathcal{B}_t\}_{t\geq 0}$ is the standard $\left(C\left([0, 1]\times[0, 1]\right)\right)^\prime$-valued Brownian motion such that $\{\mathcal{B}_t(H)\}_{t\geq 0}$ is a real-valued Brownian motion with $\mathcal{B}_0(H)=0$, ${\rm Var}\left(\mathcal{B}_t(H)\right)=t\int_0^1\int_0^1H^2(u,v)dudv$ for all $H\in C\left([0, 1]\times[0, 1]\right), t\geq 0$ and $\mathcal{U}_t$ is the linear operator from $C[0, 1]$ to $C\left([0, 1]\times [0, 1]\right)$ such that
\[
\mathcal{U}_tf(u,v)=\sqrt{\lambda(u,v)\left(\rho(t,u)+\rho(t,v)-2\rho(t,u)\rho(t,v)\right)}\left(f(u)-f(v)\right)
\]
for all $f\in C[0, 1], u,v\in [0, 1], t\geq 0$.

\quad

\textbf{Example 3} \emph{$N$-urn binary contact path process.} In this case, $a_1=a_2=a_4=1$ and $a_3=c=0$. Hence, $P_1$ is given by
\[
P_1f(u)=-b(u)f(u)+\int_0^1\lambda(v,u)f(v)dv
\]
for any $u\in [0, 1]$ and $f\in C[0, 1]$. Furthermore, by direct calculation, $P_2=P_1$ and $l_1$ is given by
\[
l_1(f)=2\int_0^1\int_0^1\rho(t,u)\rho(t,v)\lambda(u,v)f(u)dudv
\]
for any $f\in C[0, 1]$. According to the proof of Theorem \ref{theorem 2.6 fluctuation of the NurnLinearSyetem}, $\vartheta$ is given by
\[
\begin{cases}
&\frac{d}{dt}\vartheta(t,\cdot)=-b(\cdot)\vartheta(t, \cdot)+\int_0^1\vartheta(t, v)\lambda(\cdot, v)dv
+2\rho(t, \cdot)\int_0^1 \rho(t, v)\lambda(\cdot, v)dv,\\
&\vartheta_0=\phi.
\end{cases}
\]
For any $0\leq s\leq T$ and $f\in C[0, 1]$,
\[
\left[f, f\right]_s=\int_0^1f^2(u)\left[\vartheta(s,u)b(u)+\int_0^1\vartheta(s,v)\lambda(u, v)dv\right]du.
\]
As a result, by Theorem \ref{theorem 2.2 hydrodynamic limit},
\[
\lim_{N\rightarrow+\infty}\mu_t^N(f)=\int_0^1f(u)\rho(t,u)du
\]
in $L^2$ for any $f\in C[0, 1]$ and $t\geq 0$, where $\rho$ is given by
\[
\begin{cases}
&\frac{d}{dt}\rho(t, \cdot)=-b(\cdot)\rho(t, \cdot)+\int_0^1\lambda(\cdot, v)\rho(t, v)dv,\\
&\rho(0, \cdot)=\phi.
\end{cases}
\]
By Theorem \ref{theorem 2.6 fluctuation of the NurnLinearSyetem}, $V^N$ converges weakly to $V$ as $N\rightarrow+\infty$, where $V$ is the solution to
\[
\begin{cases}
&dV_t=P_1^{*}V_tdt+J_t^*d\mathcal{W}_t\text{~for~}0\leq t\leq T,\\
&V_0(f) \text{~follows~}\mathbb{N}\left(0, \int_0^1f^2(u)\phi(u)(1-\phi(u))du\right)\text{~for any~}f\in C[0, 1],\\
&V_0 \text{~is independent of~}\{C_t\}_{0\leq t\leq T},
\end{cases}
\]
where $\{\mathcal{W}_t\}_{t\geq 0}$ is the standard $\left(C[0, 1]\right)^\prime$-valued Brownian motion as in Example 1 and $J_t$ is the linear operator from $C[0, 1]$ to $C[0, 1]$ such that
\[
J_tf(u)=\sqrt{\vartheta(t,u)b(u)+\int_0^1\vartheta(t,v)\lambda(u, v)dv}f(u)
\]
for all $f\in C[0, 1], u\in [0, 1], t\geq 0$.

\appendix{}
\section{Appendix}
\subsection{The proof of Lemma \ref{lemma 2.1 ODE of Hydrodynamic}}\label{appendix A.1}
In this subsection we prove Lemma \ref{lemma 2.1 ODE of Hydrodynamic}.

\proof[Proof of Lemma \ref{lemma 2.1 ODE of Hydrodynamic}]

As we have shown in Section \ref{section two}, for any $\nu_1, \nu_2\in \left(C[0, 1]\right)^{\prime}$,
\[
\|P_1^*\nu_1-P_1^*\nu_2\|\leq \left(\|b\|_\infty\left(\|c\|_\infty+1\right)+\|\lambda\|_\infty\left(2+\sum_{i=1}^4\|a_i\|_\infty\right)\right)\|\nu_1-\nu_2\|,
\]
i.e., Equation \eqref{equ definition of hydro ODE} satisfies Lipschitz's condition. Therefore, according to Theorem 19.1.2 of \cite{Lang}, the unique solution $\mu$ to Equation \eqref{equ definition of hydro ODE} is given by
\[
\mu_t=e^{tP_1^*}\mu_0=\sum_{k=0}^{+\infty}\frac{t^k\left(P_1^*\right)^k\mu_0}{k!},
\]
where $\mu_0(du)=\phi(u)du$. According to a similar analysis, the unique solution to Equation \eqref{equ definition of hydro density} is given by
\[
\rho_t=e^{tP_8}\phi=\sum_{k=0}^{+\infty}\frac{t^kP_8^k\phi}{k!},
\]
where $P_8$ is the linear operator from $C[0, 1]$ to $C[0, 1]$ such that
\begin{align*}
P_8f(u)=&b(u)\left(c(u)-1\right)f(u)+f(u)\int_0^1\lambda(v,u)\left(a_4(v,u)-1\right)dv\\
&+f(u)\int_0^1\lambda(u,v)\left(a_1(u,v)-1\right)dv+\int_0^1\lambda(u,v)a_2(u,v)f(v)dv\\
&+\int_0^1\lambda(v,u)a_3(v,u)f(v)dv
\end{align*}
for any $f\in C[0, 1]$ and $u\in [0, 1]$. Let $\mu_{1,t}=\mu_t=e^{tP_1^*}\mu_0$ and $\mu_{2,t}(du)=\rho(t,u)du$, then for $i=1,2$,
\[
\mu_{i,t}(f)=\int_0^1\phi(u)f(u)du+\int_0^t\mu_{i,s}(P_1f)ds
\]
for any $f\in C[0, 1]$ according to the definition of $P_8$. As a result,
\[
\|\mu_{1,t}-\mu_{2,t}\|\leq C_1\int_0^t\|\mu_{1,s}-\mu_{2,s}\|ds
\]
for any $t\geq 0$, where
\[
C_1=\|b\|_\infty\left(\|c\|_\infty+1\right)+\|\lambda\|_\infty\left(2+\sum_{i=1}^4\|a_i\|_\infty\right).
\]
Therefore, by Gronwall's inequality,
\[
\|\mu_{1,t}-\mu_{2,t}\|\leq 0e^{c_1t}=0
\]
and hence $\mu_t(du)=\mu_{2,t}(du)=\rho(t,u)du$.

\qed

\subsection{Expressions of $\mathcal{M}_N$ and $\hat{\mathcal{M}}_N$}\label{appendix A.2}
In this subsection we give $\mathcal{M}_N$ and $\hat{\mathcal{M}}_N$. By \eqref{equ 3.2}, for $i\neq j$,
\begin{align*}
&\mathcal{M}_N\left((i,j), (l,k)\right)=\\
&\begin{cases}
&\frac{1}{N}\lambda\left(\frac{j}{N}, \frac{k}{N}\right)a_2\left(\frac{j}{N}, \frac{k}{N}\right)
+\frac{1}{N}\lambda\left(\frac{k}{N}, \frac{j}{N}\right)a_3\left(\frac{k}{N}, \frac{j}{N}\right)
\text{~if~}l=i\text{~and~}k\neq j,\\
&\frac{1}{N}\lambda\left(\frac{i}{N}, \frac{l}{N}\right)a_2\left(\frac{i}{N}, \frac{l}{N}\right)
+\frac{1}{N}\lambda\left(\frac{l}{N}, \frac{i}{N}\right)a_3\left(\frac{l}{N}, \frac{i}{N}\right) \text{~if~}k=j\text{~and~}l\neq i,\\
&b\left(\frac{j}{N}\right)\left(c\left(\frac{j}{N}\right)-1\right)+\frac{1}{N}\sum_{r\neq j}\lambda\left(\frac{r}{N}, \frac{j}{N}\right)
\left(a_4\left(\frac{r}{N}, \frac{j}{N}\right)-1\right)\\
&\text{\quad\quad}+\frac{1}{N}\sum_{r\neq j}\lambda\left(\frac{j}{N}, \frac{r}{N}\right)\left(a_1\left(\frac{j}{N}, \frac{r}{N}\right)-1\right)\\
&\text{\quad\quad}+\frac{1}{N}\sum_{r\neq i}\lambda\left(\frac{r}{N}, \frac{i}{N}\right)\left(a_4\left(\frac{r}{N}, \frac{i}{N}\right)-1\right)
+b\left(\frac{i}{N}\right)\left(c\left(\frac{i}{N}\right)-1\right)\\
&\text{\quad\quad}+\frac{1}{N}\sum_{r\neq i}\lambda\left(\frac{i}{N}, \frac{r}{N}\right)\left(a_1\left(\frac{i}{N}, \frac{r}{N}\right)-1\right)
\text{~if~}(l, k)=(i, j),\\
&0\text{~else}.
\end{cases}
\end{align*}
For $i=j$,
\begin{align*}
&\mathcal{M}_N\left((i,i), (l,k)\right)=\\
&\begin{cases}
&\frac{2}{N}\lambda\left(\frac{i}{N}, \frac{k}{N}\right)a_2\left(\frac{i}{N}, \frac{k}{N}\right)
+\frac{2}{N}\lambda\left(\frac{k}{N}, \frac{i}{N}\right)a_3\left(\frac{k}{N}, \frac{i}{N}\right) \text{~if~}l=i\text{~and~}k\neq i,\\
&2b\left(\frac{i}{N}\right)\left(c\left(\frac{i}{N}-1\right)\right)+\frac{2}{N}\sum_{r\neq i}\lambda\left(\frac{r}{N}, \frac{i}{N}\right)
\left(a_4\left(\frac{r}{N}, \frac{i}{N}\right)-1\right) \\
&\text{\quad\quad}+\frac{2}{N}\sum_{r\neq i}\lambda\left(\frac{i}{N}, \frac{r}{N}\right)\left(a_1\left(\frac{i}{N}, \frac{r}{N}\right)-1\right)
\text{~if~}(l, k)=(i,i),\\
&0 \text{~else}.
\end{cases}
\end{align*}
By \eqref{equ linear system Chapman-Kolmogorov}, for $i\neq j$,
\begin{align*}
&\hat{\mathcal{M}}_N\left((i,j), (l,k)\right)=\\
&\begin{cases}
&\mathcal{M}_N\left((i,j), (l,k)\right) \text{~if~}(l,k)\neq (i,j), (i,i), (j,j),\\
&\mathcal{M}_N\left((i, j), (i, j)\right)-\frac{1}{N}\lambda\left(\frac{i}{N}, \frac{j}{N}\right)\left(a_4\left(\frac{i}{N}, \frac{j}{N}\right)-1\right)\\
&\text{\quad\quad}-\frac{1}{N}\lambda\left(\frac{j}{N}, \frac{i}{N}\right)\left(a_1\left(\frac{j}{N}, \frac{i}{N}\right)-1\right)
-\frac{1}{N}\lambda\left(\frac{j}{N}, \frac{i}{N}\right)\left(a_4\left(\frac{j}{N}, \frac{i}{N}\right)-1\right)\\
&\text{\quad\quad}-\frac{1}{N}\lambda\left(\frac{i}{N}, \frac{j}{N}\right)\left(a_1\left(\frac{i}{N}, \frac{j}{N}\right)-1\right)\\
&\text{\quad\quad}+\frac{\lambda\left(\frac{i}{N}, \frac{j}{N}\right)}{N}
\left(a_1\left(\frac{i}{N}, \frac{j}{N}\right)a_4\left(\frac{i}{N}, \frac{j}{N}\right)+a_2\left(\frac{i}{N}, \frac{j}{N}\right)a_3\left(\frac{i}{N}, \frac{j}{N}\right)-1\right)\\
&\text{\quad\quad}+\frac{\lambda\left(\frac{j}{N}, \frac{i}{N}\right)}{N}
\left(a_1\left(\frac{j}{N}, \frac{i}{N}\right)a_4\left(\frac{j}{N}, \frac{i}{N}\right)+a_2\left(\frac{j}{N}, \frac{i}{N}\right)a_3\left(\frac{j}{N}, \frac{i}{N}\right)-1\right)\\
& \text{\quad\quad\quad\quad if~}(l,k)=(i,j),\\
& \frac{\lambda\left(\frac{i}{N}, \frac{j}{N}\right)}{N}a_1\left(\frac{i}{N}, \frac{j}{N}\right)a_3\left(\frac{i}{N}, \frac{j}{N}\right)
+\frac{\lambda\left(\frac{j}{N}, \frac{i}{N}\right)}{N}a_2\left(\frac{j}{N}, \frac{i}{N}\right)a_4\left(\frac{j}{N}, \frac{i}{N}\right)\text{~if~}
(l,k)=(i,i), \\
& \frac{\lambda\left(\frac{j}{N}, \frac{i}{N}\right)}{N}a_1\left(\frac{j}{N}, \frac{i}{N}\right)a_3\left(\frac{j}{N}, \frac{i}{N}\right)
+\frac{\lambda\left(\frac{i}{N}, \frac{j}{N}\right)}{N}a_2\left(\frac{i}{N}, \frac{j}{N}\right)a_4\left(\frac{i}{N}, \frac{j}{N}\right)\text{~if~}
(l,k)=(j,j).
\end{cases}
\end{align*}
For $i=j$,
\begin{align*}
&\hat{\mathcal{M}}_N\left((i,i), (l,k)\right)=\\
&\begin{cases}
&b\left(\frac{i}{N}\right)\left(c^2\left(\frac{i}{N}\right)-1\right)+\sum_{k\neq i}\frac{\lambda\left(\frac{i}{N}, \frac{k}{N}\right)}{N}
\left(a_1^2\left(\frac{i}{N}, \frac{k}{N}\right)-1\right)\\
&\text{\quad\quad}+\sum_{k\neq i}\frac{\lambda\left(\frac{k}{N}, \frac{i}{N}\right)}{N}\left(a_4^2\left(\frac{k}{N}, \frac{i}{N}\right)-1\right)
\text{~if~}(l,k)=(i,i),\\
&\frac{2\lambda\left(\frac{i}{N}, \frac{j}{N}\right)}{N}a_1\left(\frac{i}{N}, \frac{j}{N}\right)a_2\left(\frac{i}{N}, \frac{j}{N}\right)
+\frac{2\lambda\left(\frac{j}{N}, \frac{i}{N}\right)}{N}a_3\left(\frac{j}{N}, \frac{i}{N}\right)a_4\left(\frac{j}{N}, \frac{i}{N}\right)\\
&\text{\quad\quad\quad}\text{~if~}(l,k)=(i, j)\text{~for~}j\neq i,\\
&\frac{\lambda\left(\frac{i}{N}, \frac{j}{N}\right)}{N}a_2^2\left(\frac{i}{N}, \frac{j}{N}\right)+\frac{\lambda\left(\frac{j}{N}, \frac{i}{N}\right)}{N}a_3^2\left(\frac{j}{N}, \frac{i}{N}\right) \text{~if~}(l,k)=(j,j)\text{~for~}j\neq i.
\end{cases}
\end{align*}

\quad

\textbf{Acknowledgments.} The author is grateful to the financial support from the National Natural Science Foundation of China with grant number 11501542.

{}
\end{document}